\documentclass[times]{nmeauth}

\usepackage{moreverb}

\usepackage[colorlinks,bookmarksopen,bookmarksnumbered,citecolor=red,urlcolor=red]{hyperref}

\usepackage{bm}
\usepackage{siunitx}

\newcommand{\KRO}[1]{\textcolor{black}{#1}}
\newcommand{\MIG}[1]{\textcolor{black}{#1}}

\newcommand\BibTeX{{\rmfamily B\kern-.05em \textsc{i\kern-.025em b}\kern-.08em
T\kern-.1667em\lower.7ex\hbox{E}\kern-.125emX}}

\def\volumeyear{2017}

\begin{document}

\runningheads{K.~Olesen et al.}{A Higher Order Equilibrium Finite Element Method}

\title{A Higher Order Equilibrium Finite Element Method}

\author{K.~Olesen\affil{1}, B.~Gervang\affil{2}, J.~N.~Reddy\affil{3}, M.~Gerritsma\affil{4}}

\address{\affilnum{1} Denmark, e-mail: squidy82@hotmail.com. \break \affilnum{2}Aarhus University, Department of Engineering, Inge Lehmanns Gade 10, 8000 Aarhus, Denmark, e-mail: bge@ase.au.dk. \break \affilnum{3}Department of Mechanical Engineering, Texas A \& M University, College Station, TX 77843-3123, USA, e-mail: jnreddy@tamu.edu. \break \affilnum{4}Delft University of Technology, Faculty of Aerospace Engineering, Kluyverweg 2, 2629 HT Delft, The Netherlands, e-mail: m.i.gerritsma@tudelft.nl.}

\begin{abstract}
\KRO{In this paper a mixed spectral element formulation is presented for planar, linear elasticity. The degrees of freedom for the stress are integrated traction components, i.e. surface force components. As a result the tractions between elements are continuous.
The formulation is based on minimization of the complementary energy subject to the constraints that the stress field should satisfy equilibrium of forces and moments. The Lagrange multiplier which enforces equilibrium of forces is the displacement field and the Lagrange multiplier which enforces equilibrium of moments is the rotation. The formulation satisfies equilibrium of forces pointwise if the body forces are piecewise polynomial. Equilibrium of moments is weakly satisfied.
Results of the method are given on orthogonal and curvilinear domains and an example with a point singularity is given.}
\end{abstract}


\maketitle

\section{Introduction} \label{sec:Introduction}
The finite element method (FEM) is a common tool in the engineering and scientific community to simulate physical problems, and \cite{Zienkiewicz_1972,Oden_1976,Bathe_1976, Reddy_2006} describe the method in depth. The FEM was first developed for structural problems, but later also adapted for heat and fluid problems \cite{Reddy_2010}, and now it is applicable for many physical problems such as electromagnetism \cite{Monk_2003}. In the FEM the unknown quantities of a problem are approximated by polynomial nodal expansions defined in elements of the computational domain. This means that the equilibrium of forces is satisfied in the computational nodes \cite{Reddy_2015}, however, the governing differential equations of force equilibrium are not satisfied within an element and across inter-element boundaries pointwise \cite{Almeida_Continuity} or \cite[p.119]{Cook_2001}. Newton's third law, which states that an action should have an equal and opposite reaction, is therefore violated at such interfaces. The equations describing structural problems and its variables are very geometrical \cite{Marsden_1983}, and a more physically correct approach would be to consider the geometric aspects of the problem, and expand variables associated to geometrical structures such as lines, surfaces and volumes in addition to points, see \cite{Almeida_Continuity} or \cite[p.2795]{Almeida_1996}.

The displacement field is typically only $C^0$ across neighbouring elements, which means that the stress field is discontinuous over the common boundary.
To make it continuous the values are often averaged, but this can give erroneous results if, for instance, the elements have different material parameters. The approach in this paper expands the stresses based on the forces on the surfaces of the element, and therefore produces continuous stress fields over the element boundaries.

Local considerations of the force equilibrium were put forward by Fraeijs De Veubeke in the 1960's through dual analysis, \cite{Veubeke_1964,Veubeke_1965,Veubeke_1980}, where two simultaneous analyses are performed, a kinematic and a static admissible model, which are coupled through complementary energy principles. The static admissible model considers the tractions on the element boundaries and connections to the stresses are established. Since the stress field is coupled to the displacement field the compatibility equations are in general not satisfied, and spurious kinematic modes are generally present. These are eliminated by applying for instance a stress potential or by direct approximations of the stress field in the elements. These elements were revisited in the 1990's by Moitinho de Almeida and Freitas and a family of hybrid finite elements was developed, \cite{Almeida_1996,Almeida_1991}. Here the force equilibrium equations were enforced through self-equilibrating shape functions applied to the stress field, i.e. the divergence of the approximated stress field is zero. In addition, the body forces are accounted for through a particular solution. \MIG{To have interelement equilibrium a boundary pressure was obtained as a projection of the approximated stresses, which was then equilibrated through a weighted residual method.}

Debongnie also worked with dual analysis and in \cite{Debongnie_1995} a re-examination of the principles was presented. Recently, Santos and Moitinho de Almeida, \cite{Santos_2014}, make use of Piola-Kirchhoff projections to account for the deformations of the elements, and this approach is also pursued in the current paper. In \cite{Kempeneers_2010} tetrahedral equilibrium elements of polynomial degree one and two are developed. Through the super element procedure the spurious kinematic modes are eliminated. In \cite{Wang_2014} tractions are used to construct an admissible stress field, and the interelement equilibrium is therefore automatically ensured. A similar approach is used in the present paper, however, the element surface forces are used as geometric degrees of freedoms (DOFs).

\KRO{In the current paper continuity of the tractions between elements is imposed strongly in a mixed finite element formulation. An alternative approach is to reconstruct the tractions in a post-processing step such that continuity of the inter-element tractions is restored. Such techniques may be found in \cite[Ch.8]{BookAlmeidaMaunder}, \cite[Ch.6]{AinsworthOden} and in \cite{LadevezeMaunder,ParesSantosDiez,LadevezeLeguillon}. These methods rely on a global solution and elementwise local solutions to establish continuity of the tractions \emph{a posteriori}. In the current paper this is established in \emph{a priori} and one global solve is needed.}

In \cite{BochevHyman,BonelleErn,Tonti,Kreeft_2011,Kreeft_stokes} it was shown that conservation laws can be exactly represented at the discrete level, while the numerical approximation occurs in the constitutive equations. These references consider scalar conservation equations like conservation of mass. In this paper, we extend these ideas to vector-valued conservation laws, such as conservation of force equilibrium in continuum mechanics. The method is based on the \emph{mimetic spectral element method} presented in \cite{Kreeft_2011,Kreeft_stokes}, where the variables of the problem are considered as real valued differential $k$-forms, which we associate to geometrical objects of dimension $k$. The geometry of elastic problems in the discrete setting has been considered in \cite{Douglas_2006}, where the elasticity complex was defined. In the work of Yavari, \cite{Yavari_2008}, the variables of an elastic problem are represented as \emph{vector- and covector-valued-differential $k$-forms}, which means that the variables are associated to geometrical objects of dimension $k$, which maps into linear vector spaces. In Yavari's work the constitutive relations expressed in vector-valued exterior calculus are more complicated than the real-valued cases. For further developments of Yavari's work see \cite{Angoshtari_2013,Yavari_2013,Angoshtari_2014}, where the last reference discusses the elastic complex in detail. In the conventional structural FEM there is also a coupling to the geometry, as shown by Reddy and Srinivasa \cite{Reddy_2015}, where it was derived that the edges of the elements act like trusses.

The method introduced in this paper does not require self-equilibrating shape functions, and the inter-element equilibrium is naturally satisfied. We introduce integral values of the stress field, i.e. surface forces on the boundaries of the element, to expand the stresses, and thereby enforce translational force equilibrium both internally and between elements. We expand the body force density field using the volume integral values of the projected body force density field, so when body forces are present we only satisfy the force equilibrium equations in a finite volume setting. However, when no body forces are stated or piecewise polynomial body forces are given, the translational force equilibrium is exactly satisfied.

The outline of this paper is as follows. In Section~\ref{sec:Gov_equations} the governing equations of a linear elastic problem are given. In Section~\ref{sec:min_problem} the problem is presented as a constrained minimization problem, while the equilibrium of forces is explained in a discrete sense in Section~\ref{sec:eq_force}. In Section~\ref{sec:expan} the expansion polynomials are presented and applied to the minimization equations. In Section~\ref{sec:discussion} we briefly summarize the steps taken in the paper so far. In Section~\ref{sec:Results_1} results obtained on an orthogonal grid are presented. The method is extended to curvilinear elements in Section~\ref{sec:transformations} and results for non-orthogonal grids are presented in Section~\ref{sec:Results_2}.  In Section~\ref{sec:Comp_energy} it is shown that the complementary energy converges from above in line with the findings done by Fraeijs~de~Veubeke \cite{Veubeke_1964,Veubeke_1965} and the equilibrium of moments is discussed. \KRO{The method is compared to a traditional displacement based FE method in Section~\ref{sec:FE_comp} and concluding remarks are given in Section~\ref{sec:conclusions}}.

\section{The governing equations} \label{sec:Gov_equations}
Consider a bounded domain $\Omega$ with boundary $\Gamma = \Gamma_u \cup \Gamma_t$, where $\Gamma_u \cap \Gamma_t = \emptyset$. Along $\Gamma_u$ the displacements, $\bm{u}$, are prescribed, whereas along $\Gamma_t$ the surface tractions, $\bm{t}$ are given. The Cartesian components of the displacement field are denoted by $\bm{u} = (u_1,u_2,u_3)$. The stress field will be denoted by $\bm{\sigma}$ and the strain field by $\bm{\varepsilon}$. The governing equations, written in Einstein notation, in elastostatics may be divided into three different equations:
\begin{itemize}
\item The compatibility equations
\begin{subequations} \label{eq:compatibility equations}
\begin{align}
\varepsilon_{ij} = \frac{1}{2} \left( \frac{\partial}{\partial x_j} u_i + \frac{\partial}{\partial x_i} u_j \right) \quad \mbox{ in } \Omega& \;, \label{eq:strain}\\
\omega_{ij} = \frac{1}{2} \left( \frac{\partial}{\partial x_i} u_j - \frac{\partial}{\partial x_j} u_i \right) \quad \mbox{ in } \Omega& \;,\label{eq:rotation}\\
u_i = \bar{u}_i \quad \mbox{ along } \Gamma_u& \;,
\end{align}
\end{subequations}
where $\omega_{ij}$ are the components of the field of rotations and $\bar{u}_i$ are the known displacement components along $\Gamma_u$.
\item The equilibrium equations
\begin{subequations} \label{eq:equilibrium equations}
\begin{align}
\frac{\partial}{\partial x_i} \sigma_{ij} + f_j = 0 \quad \mbox{ in } \Omega& \;, \label{eq:cons_lin_mom}\\
\sigma_{ij} - \sigma_{ji} = 0 \quad \mbox{ in } \Omega& \;, \label{eq:cons_ang_mom}\\
\sigma_{ij} \cdot n_i = \bar{t}_j \quad \mbox{ along } \Gamma_t& \;,\label{eq:traction}
\end{align}
\end{subequations}
where $f_j$ are the body force components in $\Omega$, $n_i$ are the unit normal vector components along $\Gamma_t$ and $\bar{t}_j$ are the known traction components along $\Gamma_t$.
\item The constitutive equations
\end{itemize}
\begin{equation}
\varepsilon_{ij} = C_{ijkl} \, \sigma_{kl} \;, \label{eq:stress_strain_rel_strain_exp}
\end{equation}
where $C_{ijkl}$ is the compliance tensor.

This system of equations differs slightly from the equations generally encountered for linear elasticity; first of all we work with six shear stress components instead of three and explicitly add equilibrium of moments, \eqref{eq:cons_ang_mom}. Furthermore, we introduce rotation, $\bm{\omega}$, \cite[(2.1)]{Veubeke_1964} and \cite[\S 83]{Timoshenko_1982}, as additional DOFs. The reason is that we treat equilibrium of forces,
\eqref{eq:cons_lin_mom} and equilibrium of moments, \eqref{eq:cons_ang_mom} as constraints in the minimization of the
complementary energy. The Lagrange multipliers which will enforce these constraints are the displacement components $u_i$ for \KRO{equilibrium of forces} and the rotation components, $\omega_{ij}$, for \KRO{equilibrium of moments}.

\section{The minimization of the complementary energy} \label{sec:min_problem}
In this section we will formulate the Lagrangian of the minimization problem, but first the governing equations are formulated in engineering notation. From \eqref{eq:strain} and \eqref{eq:rotation} it is clear that
\begin{equation*}
\varepsilon_{ij} = \frac{\partial}{\partial x_i} u_j - \omega_{ij} = \frac{\partial}{\partial x_i} u_j + \omega_{ji} \;,
\end{equation*}
so the strain vector can be written as
\begin{equation}
\bm{\varepsilon} = \bm{D}^T \bm{u} + \bm{R}^T \bm{\omega} \;, \label{eq:strain_mat}
\end{equation}
\begin{equation*}
\bm{\varepsilon} = \begin{Bmatrix}
\varepsilon_{11} & \varepsilon_{21} & \varepsilon_{31} & \varepsilon_{12} & \varepsilon_{22} & \varepsilon_{32} & \varepsilon_{13} & \varepsilon_{23} & \varepsilon_{33}
\end{Bmatrix}^T \;,
\end{equation*}
and
\begin{equation}
\bm{D}^T = \begin{bmatrix}
\frac{\partial}{\partial x_1} & \frac{\partial}{\partial x_2} & \frac{\partial}{\partial x_3} & 0 & 0 & 0 & 0 & 0 & 0 \\
0 & 0 & 0 & \frac{\partial}{\partial x_1} & \frac{\partial}{\partial x_2} & \frac{\partial}{\partial x_3} & 0 & 0 & 0 \\
0 & 0 & 0 & 0 & 0 & 0 & \frac{\partial}{\partial x_1} & \frac{\partial}{\partial x_2} & \frac{\partial}{\partial x_3}
\end{bmatrix}^T \;, \label{eq:grad_mat}
\end{equation}
being the gradient matrix, where $T$ denotes the transposed operator, $\bm{R}^T$ is a matrix given by
\begin{equation}
\bm{R}^T = \begin{bmatrix}
0 & 0 & 0 & 0 & 0 & 1 & 0 & -1 & 0 \\
0 & 0 & 1 & 0 & 0 & 0 & -1 & 0 & 0 \\
0 & 1 & 0 & -1 & 0 & 0 & 0 & 0 & 0 \\
\end{bmatrix}^T \;, \label{eq:rot_mat}
\end{equation}
and the vector with rotation components is given by
\begin{equation*}
\bm{\omega} = \begin{Bmatrix}
\omega_{1} & \omega_{2} & \omega_{3}
\end{Bmatrix}^T = \begin{Bmatrix}
\omega_{23} & \omega_{13} & \omega_{12}
\end{Bmatrix}^T \;.
\end{equation*}

Writing the stress vector as
\begin{equation*}
\bm{\sigma} = \begin{Bmatrix}
\sigma_{11} & \sigma_{21} & \sigma_{31} & \sigma_{12} & \sigma_{22} & \sigma_{32} & \sigma_{13} & \sigma_{23} & \sigma_{33}
\end{Bmatrix}^T \;,
\end{equation*}
then \eqref{eq:cons_lin_mom} can be expressed as
\begin{equation*}
\bm{D} \bm{\sigma} + \bm{f} = 0 \;,
\end{equation*}
where $\bm{D}$ is the divergence matrix, which is the adjoint of the gradient matrix in \eqref{eq:grad_mat}.

The equilibrium of moments is written as
\begin{equation*}
-\bm{R} \bm{\sigma} = 0 \;,
\end{equation*}
where $\bm{R}$ is the transpose of the matrix in \eqref{eq:rot_mat}.

The complementary energy is, according to \cite[\S 3]{Veubeke_1964} and \cite[(8)]{Kempeneers_2010}, given by
\begin{equation*}
CE(\bm{\sigma}) = \frac{1}{2} \int_{\Omega} \bm{\sigma}^T \bm{C} \bm{\sigma} \, \mathrm{d} \Omega - \int_{\Gamma_u} \bm{t}^T \bar{\bm{u}} \, \mathrm{d} \Gamma \;, \label{eq:comp_energy}
\end{equation*}
so the Lagrange functional for the minimization is
\begin{align}
\mathcal{L}(\bm{\sigma},\bm{u},\bm{\omega};\bm{f},\bar{\bm{u}}) = & \frac{1}{2} \int_{\Omega} \bm{\sigma}^T \bm{C} \bm{\sigma} \, \mathrm{d} \Omega - \int_{\Gamma_u} \bm{t}^T \bar{\bm{u}} \, \mathrm{d} \Gamma + \nonumber \\
&\int_{\Omega} \bm{u}^T \left( \bm{D} \bm{\sigma} + \bm{f} \right) \,\mathrm{d} \Omega - \int_{\Omega} \bm{\omega}^T \bm{R} \bm{\sigma} \,\mathrm{d} \Omega \;.
\label{eq:Elastic_lagrangian}
\end{align}

\vskip 0.3cm
\noindent
{\bf Theorem 1}
\textit{The stationary points of the Lagrangian \eqref{eq:Elastic_lagrangian} solve the linear elastic problem \eqref{eq:equilibrium equations} and \eqref{eq:stress_strain_rel_strain_exp}}.

\vskip 0.3cm
\noindent
{\bf Proof} 
\MIG{The proof is straightforward. Taking variations of the Lagrangian with respect to the rotation $\omega$ yields the symmetry of the stress tensor/equilibrium of moments. If we take variations of the functional with respect to the displacement field $\bm{u}$ we retrieve equilibrium of forces. Taking variations with respect to the stress tensor of the Lagrangian $\mathcal{L}$ gives the constitutive stress-strain relation.}

\vskip 0.3cm
\noindent
\textbf{Remark 1}
The use of rotation $\omega$ as a Lagrange multiplier to enforce the rotational equilibrium condition $\tau_{yx}=\tau_{xy}$ was already discussed by Fraeijs de Veubeke, \cite{Veubeke_1980,Veubeke_Millard_1976}. See also \cite{Bertoti} for this approach. The use of vorticity as an independent unknown in flow problems is addressed in \cite{MEEVC}.

\vskip 0.3cm
\noindent
\textbf{Remark 2}
The idea of using kinematic variables as Lagrange multipliers to enforce dynamic constraints is also applied in \cite{GerritsmaPhillips}.

\vskip 0.3cm
\noindent
\textbf{Remark 3}
In \cite{Veubeke_1964} the body force densities were incorporated through a known particular stress field, however, in this case it is incorporated directly through \eqref{eq:Elastic_lagrangian}.

\section{The equilibrium of forces} \label{sec:eq_force}
Through the use of the expansion polynomials presented in Section~\ref{sec:expan} the expansion coefficients of the approximated fields of $\bm{\sigma}$ and $\bm{f}$ are the surface forces and body forces acting on finite volumes. By considering the forces instead of the stresses and body force densities we will end up with algebraic relations, which only depend of the topology of the element, i.e no interpolation is needed.

\begin{figure}
\centering
\includegraphics[width=0.9\textwidth]{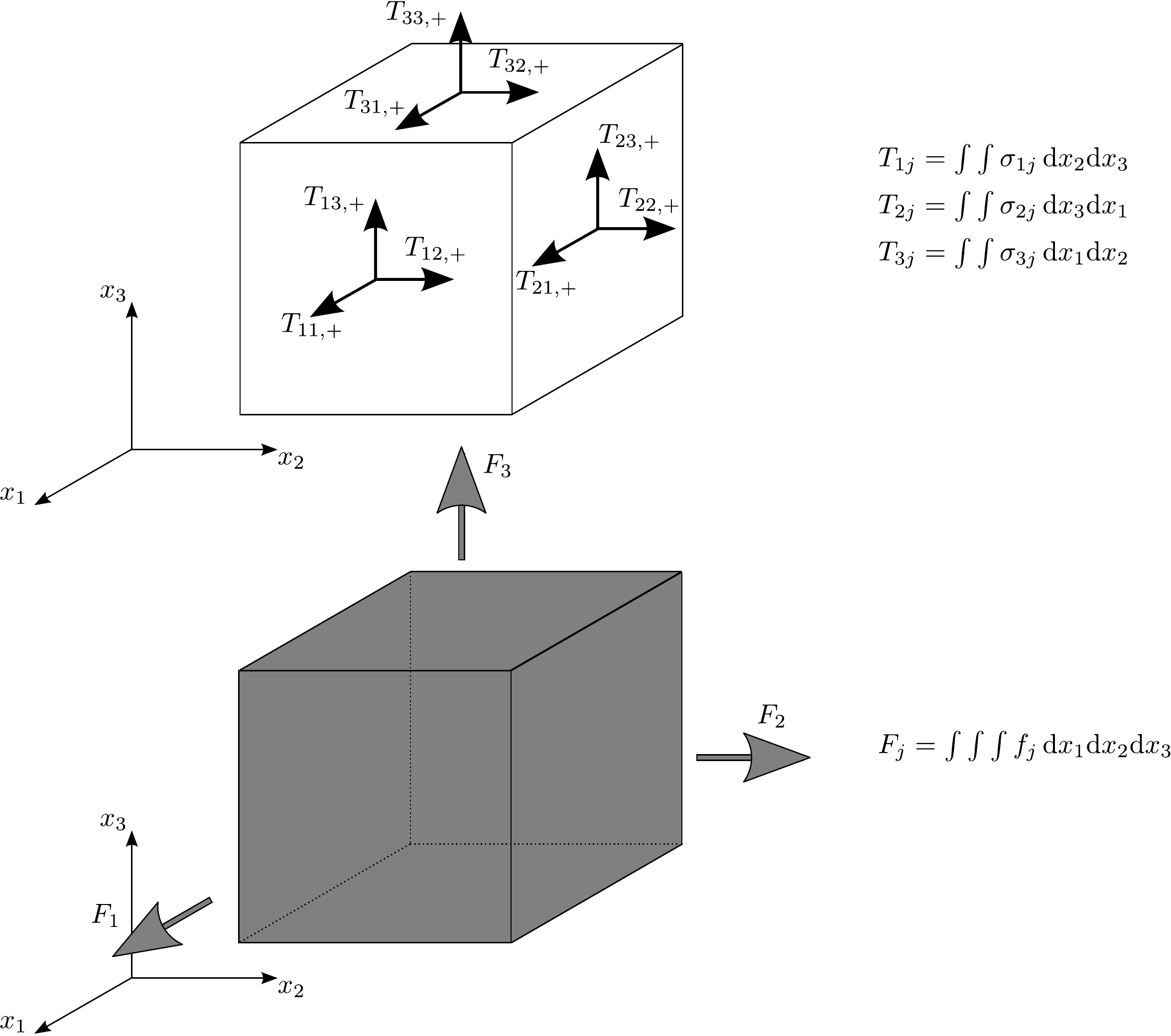}
\caption{The forces on a cube in $\mathbb{R}^3$. The top picture illustrates the components of the surface forces on the individual surfaces, $T_{i j}$. The bottom picture shows the body force components, $F_j$, acting on the volume of the element. In addition the relation $T_{i j}$ and $\sigma_{ij}$ as well as $F_j$ and $f_j$ are highlighted.}
\label{fig:Forces_nice}
\end{figure}

Consider the cube in $\mathbb{R}^3$ depicted in Figure~\ref{fig:Forces_nice}, where external forces act on the 6 boundaries of the cube and a body force acts on the volume.
The surface forces are decomposed into components in the basis directions, and are denoted by $T_{ij}$, where index $i$ denotes the global Cartesian directions of the boundary, and index $j$ denotes the direction of the surface force component. Furthermore a plus sign after the indices denotes the boundary having a unit surface normal vector in the positive basis direction, while a minus sign has the opposite meaning. The body force is also divided into the components $F_j$. Note that the integral relation between $T_{i j}$ and $\sigma_{ij}$ as well as $F_j$ and $f_j$ are also given in Figure~\ref{fig:Forces_nice}. The equilibrium of forces is written as
\begin{equation}
T_{1j,+} - T_{1j,-} + T_{2j,+} - T_{2j,-} + T_{3j,+} - T_{3j,-} + F_j = 0 \;, \label{eq:force_eq_C}
\end{equation}
for $j=1,\ldots,3$. Consider a domain $\Omega$ in $\mathbb{R}^3$ and dividing this in a number of sub-domains, i.e. elements, then the equilibrium of forces in $\Omega$ must be given by \eqref{eq:force_eq_C} for all the individual elements. If the individual force components are numbered and arranged in the column vectors $\bm{\Delta}_T$ for the surface force components and $\bm{\Delta}_F$ for body force components, then the equilibrium of forces can be written as
\begin{equation}
\bm{\mathcal{D}} \bm{\Delta}_T = - \bm{\Delta}_F \;, \label{eq:force_phys}
\end{equation}
where $\bm{\mathcal{D}}$ is a matrix only containing the numbers $-1$, $1$ and $0$ resulting from \eqref{eq:force_eq_C}.

\begin{figure}
\centering
\includegraphics[width=1.0\textwidth]{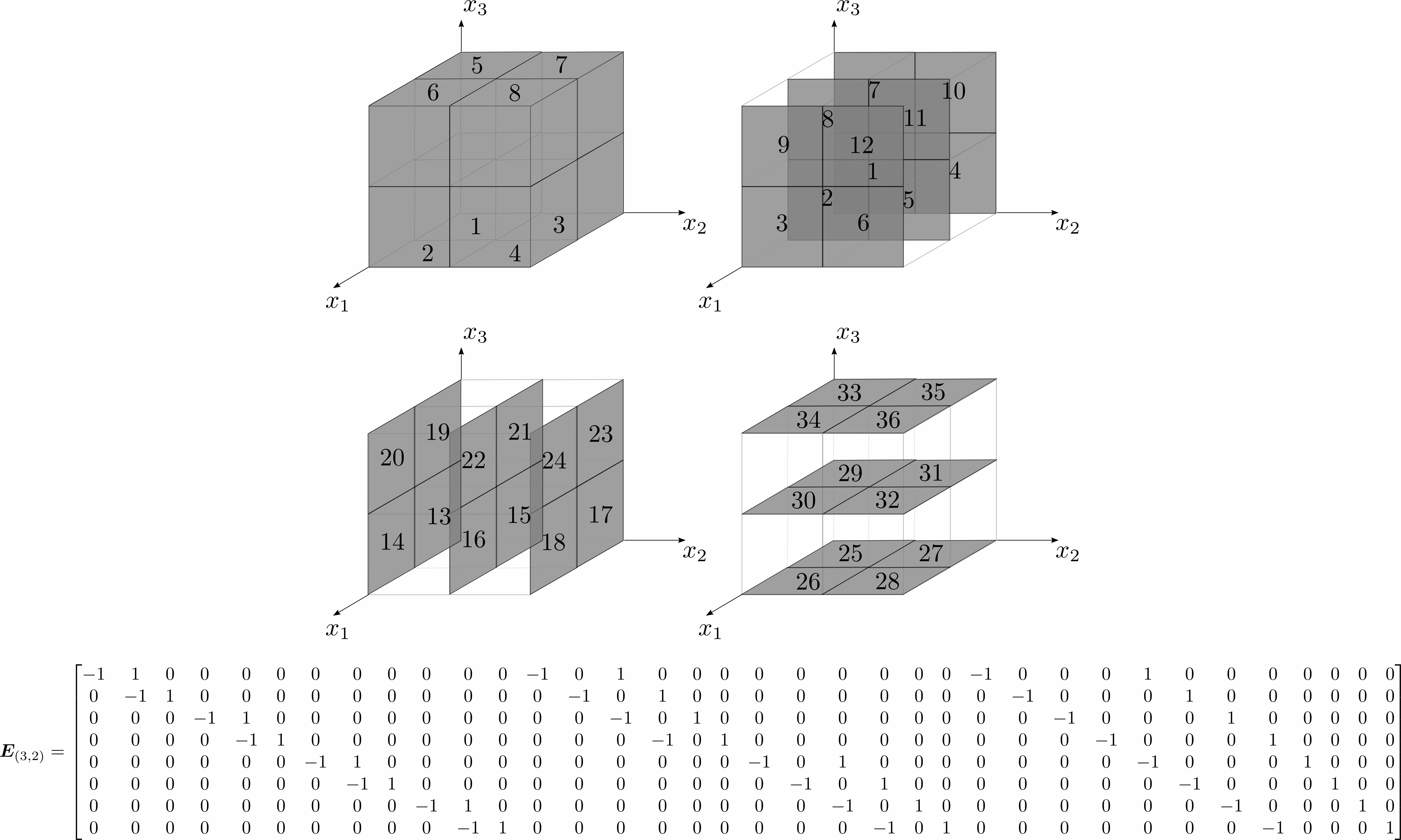}
\caption{Example of the numbering of the elements and their bounding surfaces. The numbering is increasing first in the $x_1$- then in $x_2$- and lastly in the $x_3$-direction. The illustrations show the numbering of: top left: The elements, top right: Bounding surfaces normal to $x_1$, lower left: Bounding surfaces normal to $x_2$, and lower right: Bounding surfaces normal to $x_3$.}
\label{fig:Fluxes_connection}
\end{figure}

From \eqref{eq:force_eq_C} we see that the matrix $\bm{\mathcal{D}}$ only depends on the connectivity of the elements and can be written as
\begin{equation*}
\bm{\mathcal{D}} = \begin{bmatrix}
\bm{E}_{(3,2)} & 0 & 0 \\
0 & \bm{E}_{(3,2)} & 0 \\
0 & 0 & \bm{E}_{(3,2)}
\end{bmatrix} \;, \label{eq:div_mat}
\end{equation*}
where $\bm{E}_{(3,2)}$ is used as a discrete divergence operator in \cite{Kreeft_2011,Kreeft_stokes}, and is also denoted an incidence matrix. Figure~\ref{fig:Fluxes_connection} shows a domain which is divided into 8 elements, and the numbering of these as well as its bounding surfaces are illustrated, which results in the highly sparse matrix $\bm{E}_{(3,2)}$ written out under the pictures.
The surface force components $T_{i,j\pm}$ and the integrated body force component $F_j$, see Figure~\ref{fig:Fluxes_connection}, in \eqref{eq:force_eq_C} will be explicitly used as expansion coefficients in the high order polynomial approximation to be discussed in Section~\ref{sec:expan}, see equations \eqref{eq:stress_expan} and \eqref{eq:vol_expan}.

\section{Expansion Polynomials} \label{sec:expan}
\begin{figure}
\centering
\includegraphics[width=0.5\textwidth]{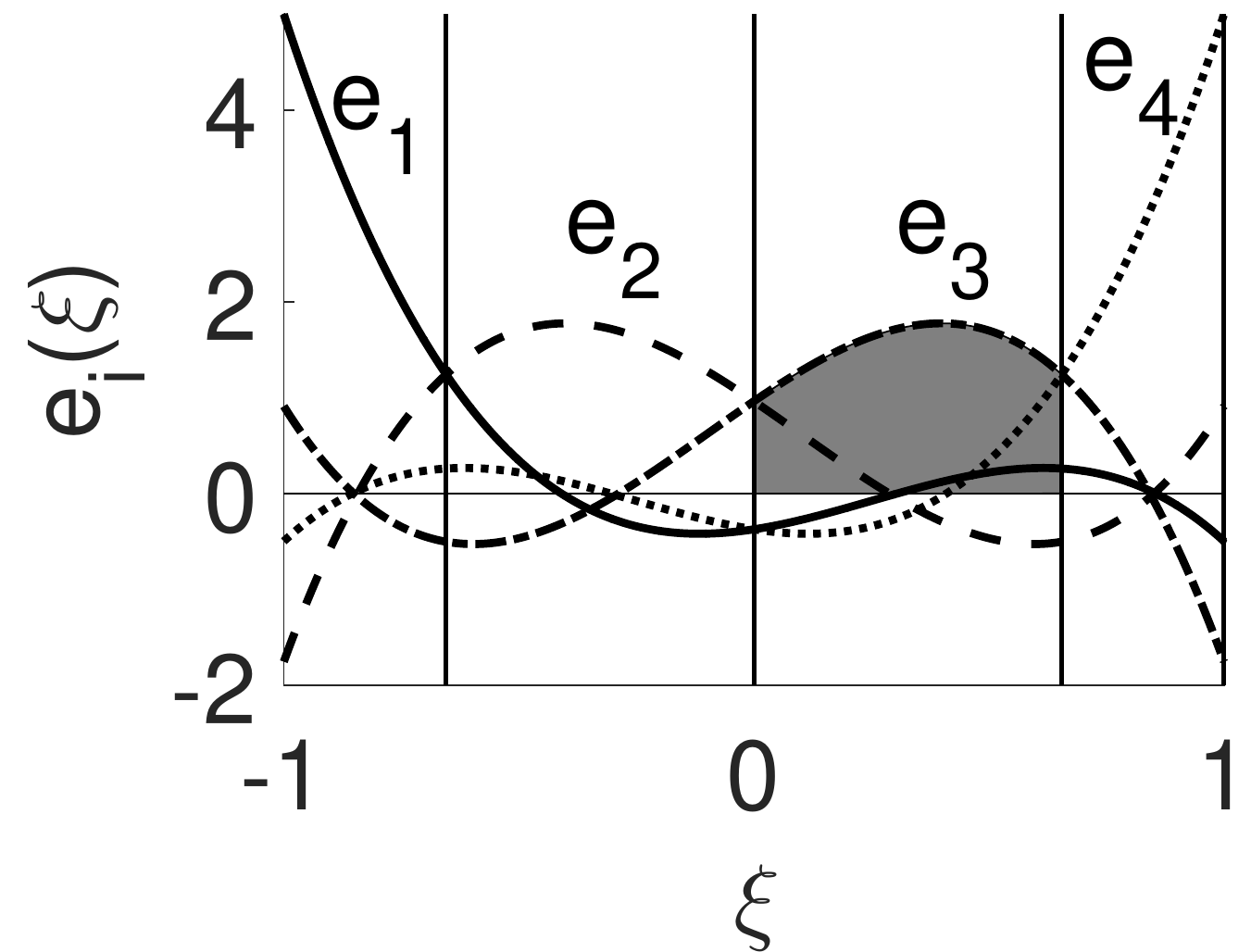}
\caption{Edge polynomials of polynomial degree $N-1$ associated to the Gauss-Lobatto-Legendre
points for $N=4$. The property $\int_{\xi_{k-1}}^{\xi_k} e_i(\xi) = \delta_{i,k}$ is highlighted for the $e_3(\xi)$ polynomial.}
\label{fig:edge_poly}
\end{figure}

This section describes the expansion polynomials and the resulting linear equation system describing the problem. From Figure~\ref{fig:Forces_nice} it is seen that a surface force component is the surface integral of the corresponding stress component. Therefore, it is evident that if the expansion of the stress components is based on integral values then we will have a connection between the discrete force values $T_{ij}$ in \eqref{eq:force_phys} and the expanded stress components in $\sigma_{ij}^h$. Such polynomials are derived in \cite{Gerritsma_2011} and are named \emph{edge polynomials}. Just as Lagrange polynomials of polynomial degree $N$, defined by
\begin{equation*}
h_i(\xi) = \frac{\prod_{j=0,j \neq i}^N(\xi-\xi_j)}{\prod_{j=0,j \neq i}^N(\xi_i-\xi_j)} \;,
\end{equation*}
with $\xi_i$, $i=0,\ldots,N$ being nodal points, have the property
\begin{equation}
h_i(\xi_j) = \delta_{i,j} \;,
\label{eq:Kronecker_h}
\end{equation}
then edge polynomials of polynomial degree $(N-1)$ are given by
\begin{equation*}
e_i(\xi) = - \sum_{k=0}^{i-1} \frac{\mathrm{d} h_k(\xi)}{\mathrm{d} \xi} \;,\;\;\; i=1,\ldots,N \;,
\end{equation*}
have the property:
\begin{equation}
\int_{\xi_{k-1}}^{\xi_k} e_i(\xi) = \delta_{i,k} = \left\lbrace \begin{matrix} 1 \quad \mbox{if } i=k \\ 0 \quad \mbox{if } i \neq k \end{matrix} \right. \;.
\label{eq:Kronecker_e}
\end{equation}
This is illustrated in Figure~\ref{fig:edge_poly}.

Let $\phi(\xi)$ in $\mathbb{R}$ be a function with the expansion
\begin{equation*}
\phi^h(\xi) = \sum\limits_{i=0}^N \phi_i h_i(\xi) \;,
\end{equation*}
where $\phi_i$ is the nodal expansion coefficients, then $\frac{d\phi^h}{d\xi}$ is according to \cite{Gerritsma_2011} given by
\begin{equation}
\frac{d\phi^h}{d\xi} = \sum_{i=0}^N \phi_i \frac{dh_i(\xi)}{d\xi} =
\sum_{i=1}^N (\phi_i - \phi_{i-1}) e_i(\xi) \;.
\label{eq:expan_diff}
\end{equation}

Just as the Lagrange polynomials expand functions based on nodal values the edge polynomials expand functions based on integral values along lines. Expansions based on surfaces can be created through tensor products between edge polynomials in two directions and expansions based on volumes can be represented through tensor products between edge polynomials in three directions. The expansion coefficients are based on integral values of the expanded stress field and body force density field, respectively, i.e. surface force components and body force components.

\subsection{Dynamic variables}
The stress components $\sigma_{nm}$ are expanded by
\begin{subequations} \label{eq:stress_expan}
\begin{align}
\sigma_{1 m}^{s,h}(\xi_1,\xi_2,\xi_3) =& \sum\limits_{i=0}^{N} \sum\limits_{j=1}^{N} \sum\limits_{k=1}^{N} \left( T_{1 m}^s \right)_{i,j,k} h_i(\xi_1) e_j(\xi_2) e_k(\xi_3) \;,
\\
\sigma_{2 m}^{s,h}(\xi_1,\xi_2,\xi_3) =& \sum\limits_{i=1}^{N} \sum\limits_{j=0}^{N} \sum\limits_{k=1}^{N} \left( T_{2 m}^s \right)_{i,j,k} e_i(\xi_1) h_j(\xi_2) e_k(\xi_3) \;,
\\
\sigma_{3 m}^{s,h}(\xi_1,\xi_2,\xi_3) =& \sum\limits_{i=1}^{N} \sum\limits_{j=1}^{N} \sum\limits_{k=0}^{N} \left(T_{3 m}^s \right)_{i,j,k} e_i(\xi_1) e_j(\xi_2) h_k(\xi_3) \;,
\end{align}
\end{subequations}
where the superscript $s$ specifies that it is valid in element $s$. If $f_m$ is the body force density component then this is expanded based on the body force component $F_m$, i.e.
\begin{equation}
f_m^{s,h}(\xi_1,\xi_2,\xi_3) = \sum\limits_{i=1}^{N} \sum\limits_{j=1}^{N} \sum\limits_{k=1}^{N} \left( F_{m}^s \right)_{i,j,k} e_i(\xi_1) e_j(\xi_2) e_k(\xi_3) \;, \label{eq:vol_expan}
\end{equation}
with
\begin{equation*}
F_{m}^s = \int\limits_{\Omega_s} f_{m} \, \mathrm{d} V \;.
\end{equation*}

Gathering the $9 N^2 (N+1)$ expansion coefficients from \eqref{eq:stress_expan} in the column vector $\bm{\Delta}_T^s$ then \eqref{eq:stress_expan} can be written as
\begin{equation*}
\bm{\sigma}^{s,h} = \bm{\Psi}_2 \bm{\Delta}_T^s \;,
\end{equation*}
where the matrix $\bm{\Psi}_2$ (this notation is inspired from \cite{Reddy_2006}) is of size $9 \times 9 N^2 (N+1)$ and contains all the expansion polynomials. Similarly \eqref{eq:vol_expan} is written as
\begin{equation*}
\bm{f}^{s,h} = \bm{\Psi}_3 \bm{\Delta}_F^s \;,
\end{equation*}
where $\bm{\Psi}_3$ is a matrix of size $3 \times 3 N^3$ and $\bm{\Delta}_F^s$ is a column vector containing the $3 N^3$ expansion coefficients of the element. Applying a global assembly known from conventional FEM, the stress components and the body force densities are described by the global system
\begin{equation*}
\bm{\sigma}^{h} = \bm{\Psi}_2^G \bm{\Delta}_T \;,
\end{equation*}
and
\begin{equation*}
\bm{f}^{h} = \bm{\Psi}_3^G \bm{\Delta}_F \;,
\end{equation*}
respectively.

The divergence of the expanded stress field is given by applying \eqref{eq:expan_diff} to \eqref{eq:stress_expan}, which will have the consequence that all the expansion polynomials will consist of edge polynomials. If the surface forces are numbered as in Figure~\ref{fig:Forces_nice}
then $\bm{D} \bm{\sigma}^h$, which will be used in \eqref{eq:approx_lin_mom_sum}, is written as
\begin{align*}
\bm{D} \bm{\sigma}^h = \bm{\Psi}_3^G \bm{\mathcal{D}} \bm{\Delta}_T \;.
\end{align*}

\subsection{Kinematic variables}
The displacement components are expanded using Lagrange polynomials, however, since the displacement components act as Lagrange multipliers to enforce \eqref{eq:cons_lin_mom} we need an equal amount of displacement DOFs as the number of equations in \eqref{eq:force_phys}. This is naturally satisfied if the displacement DOFs are located in the Gauss-Legendre (GL) points and the stress and body force components are expanded on the surfaces and volumes spanned by the Gauss-Lobatto-Legendre (GLL) points. The GLL grid, where the dynamical variables are expanded, is called the primal grid and the GL grid is called the dual grid and this is where the kinematics are described, see \cite{BochevHyman,Kreeft_2011,Bossavit,Hirani,Lipnikova_2014} for the use of dual grids. Let $\xi_i$, $i=0,\ldots,N$ be the GLL points of polynomial degree $N$ and $\tilde{\xi}_i$, $i=0,\ldots,N-1$ the GL points. Note that $\xi_i< \tilde{\xi}_i < \xi_{i+1}$, for $i=0,\ldots,N-1$. The Lagrange polynomials associated with the GLL points will be denoted by $h_i(\xi)$ and the Lagrange polynomials associated with the GL points will be denoted by $\tilde{h}_i(\xi)$. For more details see \cite{Canuto_2006}. The edge polynomial $e_i(\xi)$ is a polynomial of degree $N-1$ and $\tilde{e}_i(\xi)$ is a polynomial of degree $N-2$. The $m^{th}$ component of the displacement is expanded in each element of the dual grid as
\begin{equation}
u_m^{s,h}(\xi_1,\xi_2,\xi_3) = \sum\limits_{i=0}^{N-1} \sum\limits_{j=0}^{N-1} \sum\limits_{k=0}^{N-1} \left( u_m \right)_{i,j,k} \tilde{h}_i(\xi_1) \tilde{h}_j(\xi_2) \tilde{h}_k(\xi_3) \;, \label{eq:disp_expan}
\end{equation}
This is expressed in vector-form as
\begin{equation*}
\bm{u}^{s,h} = \tilde{\bm{\Psi}}_0 \bm{\Delta}_u^s \;,
\end{equation*}
where $\tilde{\bm{\Psi}}_0$ is a matrix of size $3 \times N^3$ with the expansion polynomials of \eqref{eq:disp_expan}, and the $N^3$ displacement DOFs are gathered in $\bm{\Delta}_u^s$.

The components of the rotation field, $\bm{\omega}$, function as Lagrange multipliers to enforce \eqref{eq:cons_ang_mom} and the components of $\bm{\omega}$ are therefore expanded on the GL grid, i.e.
\begin{subequations} \label{eq:expan_rot}
\begin{align}
\omega_{1}^{s,h}(\xi_1,\xi_2,\xi_3) =& \sum\limits_{i=0}^{N-1} \sum\limits_{j=0}^{N-1} \sum\limits_{k=0}^{N-1} \left( O_{1}^s \right)_{i,j,k} \tilde{h}_i(\xi_1) \tilde{h}_j(\xi_2) \tilde{h}_k(\xi_3) \;,
\\
\omega_{2}^{s,h}(\xi_1,\xi_2,\xi_3) =& \sum\limits_{i=0}^{N-1} \sum\limits_{j=0}^{N-1} \sum\limits_{k=0}^{N-1} \left( O_{2}^s \right)_{i,j,k} \tilde{h}_i(\xi_1) \tilde{h}_j(\xi_2) \tilde{h}_k(\xi_3) \;,
\\
\omega_{3}^{s,h}(\xi_1,\xi_2,\xi_3) =& \sum\limits_{i=0}^{N-1} \sum\limits_{j=0}^{N-1} \sum\limits_{k=0}^{N-1} \left(O_{3}^s \right)_{i,j,k} \tilde{h}_i(\xi_1) \tilde{h}_j(\xi_2) \tilde{h}_k(\xi_3) \;,
\end{align}
\end{subequations}
where
\[ \left( O_{1}^s \right)_{i,j,k} = \omega_{1}^{s,h}(\left( \tilde{\xi}_1 \right)_{i},\left( \tilde{\xi}_2 \right)_{j},\left( \tilde{\xi}_3 \right)_{k})\;,\;\;\;
\left( O_{2}^s \right)_{i,j,k} = \omega_{2}^{s,h}(\left( \tilde{\xi}_1 \right)_{i},\left( \tilde{\xi}_2 \right)_{j},\left( \tilde{\xi}_3 \right)_{k})\;,\]
\[ \left( O_{3}^s \right)_{i,j,k} = \omega_{3}^{s,h}(\left( \tilde{\xi}_1 \right)_{i},\left( \tilde{\xi}_2 \right)_{j},\left( \tilde{\xi}_3 \right)_{k}) \;.
\]
This is expressed in vector-form as
\begin{equation*}
\bm{\omega}^{s,h} = \tilde{\bm{\Psi}}_{0} \bm{\Delta}_O^s \;,
\end{equation*}
where $\tilde{\bm{\Psi}}_{0}$ is a matrix of size $3 \times N^3$ with the expansion polynomials of \eqref{eq:expan_rot}, and the $N^3$ rotation DOFs are gathered in $\bm{\Delta}_O^s$.

Let $\mathcal{L}(\bm{\sigma}^h,\bm{u}^h,\bm{\omega}^h;\bm{f}^h,\bar{\bm{u}}^h)$ be the Lagrange functional in \eqref{eq:Elastic_lagrangian} based on approximated fields, where the superscript $h$ denotes that it is a discrete formulation, then the corresponding variations of \eqref{eq:Elastic_lagrangian}, respectively yield the equations
\begin{subequations}
\begin{align}
\int_{\Omega} \left( \tilde{\bm{\omega}}^h \right)^T \left( \bm{R} \bm{\sigma}^h \right) \,\mathrm{d} \Omega = 0 \quad \quad & \forall \tilde{\bm{\omega}}^h \in L^2(\Omega) \;, \label{eq:approx_ang_mom} \\
\int_{\Omega} \left( \tilde{\bm{u}}^h \right)^T ( \bm{D} \bm{\sigma}^h + \bm{f}^h )\,\mathrm{d} \Omega = 0 \quad \quad & \forall  \tilde{\bm{u}}^h \in \left [ L^2(\Omega) \right ]^2\;, \label{eq:approx_lin_mom} \\
\int_{\Omega} \left( \tilde{\bm{\sigma}}^h \right)^T \left( \bm{C} \bm{\sigma}^h - \bm{\varepsilon}^h \right) \, \mathrm{d} \Omega = 0 \quad \quad & \forall  \tilde{\bm{\sigma}}^h \in \left [ H(div;\Omega) \right ]^2 \;. \label{eq:constitutive_approx}
\end{align}
\end{subequations}
Inserting \eqref{eq:strain_mat} into \eqref{eq:constitutive_approx} yields
\begin{equation*}
\int_{\Omega} \left( \tilde{\bm{\sigma}}^h \right)^T \left( \bm{C} \bm{\sigma}^h - \bm{D}^T \bm{u}^h - \bm{R}^T \bm{\omega}^h \right) \, \mathrm{d} \Omega = 0 \;.
\end{equation*}
Using integration by parts on the second term gives
\begin{equation}
\int_{\Omega} \left[ \left( \tilde{\bm{\sigma}}^h \right)^T \left( \bm{C} \bm{\sigma}^h - \bm{R}^T \bm{\omega}^h \right) + \bm{D} \tilde{\bm{\sigma}}^h \bm{u}^h \right] \, \mathrm{d} \Omega = \int_{\Gamma_u} \left( \tilde{\bm{t}}^h \right)^T \bar{\bm{u}} \, \mathrm{d} \Gamma \;, \label{eq:approx_constit}
\end{equation}
where it has been used that $\tilde{\bm{t}}=0$ on $\Gamma_t$.

Dividing $\Omega$ into $N_{e}$ non-overlapping elements, i.e. $\Omega = \sum\limits_{s=1}^{N_{e}} \Omega_s$ means that \eqref{eq:approx_ang_mom}, \eqref{eq:approx_lin_mom} and \eqref{eq:approx_constit} can be written as
\begin{subequations} \label{eq:approx_eq_sum}
\begin{align}
\sum\limits_{s=1}^{N_{e}} \int_{\Omega_s} \left( \tilde{\bm{\omega}}^h \right)^T \left( \bm{R} \bm{\sigma}^h \right) \,\mathrm{d} \Omega =& 0 \;, \label{eq:approx_ang_mom_sum} \\
\sum\limits_{s=1}^{N_{e}} \int_{\Omega_s} \left( \tilde{\bm{u}}^h \right)^T ( \bm{D} \bm{\sigma}^h + \bm{f}^h )\,\mathrm{d} \Omega =& 0 \;, \label{eq:approx_lin_mom_sum} \\
\sum\limits_{s=1}^{N_{e}} \int_{\Omega_s} \left[ \left( \tilde{\bm{\sigma}}^h \right)^T \left( \bm{C} \bm{\sigma}^h - \bm{R}^T \bm{\omega}^h \right) + \bm{D} \tilde{\bm{\sigma}}^h \bm{u}^h \right] \, \mathrm{d} \Omega =& \int_{\Gamma_u} \left( \tilde{\bm{t}}^h \right)^T \bar{\bm{u}} \, \mathrm{d} \Gamma \;. \label{eq:approx_constit_sum}
\end{align}
\end{subequations}

Inserting the expansions into \eqref{eq:approx_eq_sum} and evaluating the integrals by appropriate Gaussian quadratures yields the symmetric linear equation system

\begin{equation}
\begin{bmatrix}
\bm{H}^h & \left ( \bm{V}^h \bm{\mathcal{D}} \right )^T & -\left( \bm{R}^h \right)^T  \\
\bm{V}^h \bm{\mathcal{D}}  & \bm{0} & \bm{0} \\
-\bm{R}^h & \bm{0} &  \bm{0}
\end{bmatrix} \begin{Bmatrix}
\bm{\Delta}_T \\  \bm{\Delta}_u \\ \bm{\Delta}_O
\end{Bmatrix} = \begin{Bmatrix}
\bm{B}^h \bm{\Delta}_{\bar{u}} \\ - \bm{V}^h \bm{\Delta}_{\hat{F}} \\
\bm{0}
\end{Bmatrix} \;, \label{eq:lin_eq_sys_non_def}
\end{equation}
where
\begin{subequations}
\begin{align}
\bm{V}^h =& \sum\limits_{s=1}^{N_e} \sum\limits_{i=0}^{N-1} \sum\limits_{j=0}^{N-1} \sum\limits_{k=0}^{N-1} \left( \left( \bm{\Psi}_0^s(\tilde{\bm{\xi}}_{i,j,k}) \right)^T \bm{\Psi}_3^s(\tilde{\bm{\xi}}_{i,j,k}) \tilde{w}_i \tilde{w}_j \tilde{w}_k \right) \;, \\
\bm{H}^h =& \frac{2}{h_{el}}\sum\limits_{s=1}^{N_e} \sum\limits_{i=0}^{N} \sum\limits_{j=0}^{N} \sum\limits_{k=0}^{N} \left( \left( \bm{\Psi}_2^s(\bm{\xi}_{i,j,k}) \right)^T \bm{C} \bm{\Psi}_2^s(\bm{\xi}_{i,j,k}) w_i w_j w_k \right) \;, \label{eq:const_mat}\\
\bm{R}^h =& \sum\limits_{s=1}^{N_e} \sum\limits_{i=0}^{N} \sum\limits_{j=0}^{N} \sum\limits_{k=0}^{N} \left( \left( \bm{\Psi}_{0}^s(\bm{\xi}_{i,j,k}) \right)^T \bm{R} \tilde{\bm{\Psi}}_2^s(\bm{\xi}_{i,j,k}) w_i w_j w_k \right) \;, \\
\bm{B}^h =& \sum\limits_{s=1}^{N_{BCu}} \sum\limits_{i=0}^{N} \sum\limits_{j=0}^{N} \left( \left( \bm{\Psi}_2^s(\tilde{\bm{\xi}}_{i,j,k}) \right)^T \tilde{\bm{\Psi}}_0^s(\tilde{\bm{\xi}}_{i,j,k}) \tilde{w}_i \tilde{w}_j \right) \;.
\end{align}
\end{subequations}
In the above relations $\bm{\xi}_{i,j,k} = (\xi_{1,i},\xi_{2,j},\xi_{3,k})$ are the GLL points, $\tilde{\bm{\xi}}_{i,j,k} = (\tilde{\xi}_{1,i},\tilde{\xi}_{2,j},\tilde{\xi}_{3,k})$ are the GL points, $w_i$ are the integration weights for the GLL quadrature, $\tilde{w}_i$ are the integration weights for the GL quadrature, $N_{BCu}$ is the number of surface elements on $\Gamma_u$, and $\bm{\Delta}_{\bar{u}}$ is a column vector listing all the known displacement components in the GL points. The term $\frac{2}{h_{el}}$ in \eqref{eq:const_mat} is the scaling of the reference $2\times2\times2$ element assuming that the element is scaled proportionally in all basis directions to an element of size $h_{el}$.

Notice that the divergence of the stress field and the body force density in the equilibrium of forces in \eqref{eq:lin_eq_sys_non_def} are both expanded with $\bm{V}^h$, and therefore the expansion of the equilibrium of forces can be reduced to
\begin{equation*}
\bm{\mathcal{D}} \bm{\Delta}_T = - \bm{\Delta}_F \;,
\end{equation*}
which is equal to \eqref{eq:force_phys} meaning it is exact.

\section{Discussion}\label{sec:discussion}
In \eqref{eq:equilibrium equations} the equilibrium equations were introduced consisting of translational equilibrium, \eqref{eq:cons_lin_mom},
the rotational equilibrium, \eqref{eq:cons_ang_mom}, and the traction along the outer boundary, \eqref{eq:traction}. If translational
equilibrium is integrated over an arbitrary volume $\mathcal{V} \subset \Omega$, one obtains
\[ \int_{\partial \mathcal{V}} \sigma\cdot n\, \mathrm{d}S + \int_{\mathcal{V}} f\,\mathrm{d}\mathcal{V} = 0 \;.\]
If we divide the boundary $\partial \mathcal{V}$ in several sub-faces as, for instance, was done in Figure~1,
we end up with equation \eqref{eq:force_eq_C}. Therefore, from \eqref{eq:cons_lin_mom} we can obtain \eqref{eq:force_eq_C}, but in general, it is not
possible to retrieve \eqref{eq:cons_lin_mom} from the integral formulation \eqref{eq:force_eq_C}. But in a finite-dimensional setting, as
discussed in this paper, it is possible to satisfy the point-wise relation \eqref{eq:cons_lin_mom}.

If we use the stress representation as given in \eqref{eq:stress_expan} and take the divergence, we obtain
\begin{align*} \bm{D} \bm{\sigma}^{s,h}_m = \sum_{i,j,k=1}^N & \left [
\left ( T_{1m}^s \right )_{i,j,k} - \left ( T_{1m}^s \right )_{i-1,j,k} +
\left ( T_{2m}^s \right )_{i,j,k} - \left ( T_{2m}^s \right )_{i,j-1,k} + \right . \\
& \left .
\left ( T_{3m}^s \right )_{i,j,k} - \left ( T_{3m}^s \right )_{i,j,k-1} \right ]
e_i(\xi_1) e_j(\xi_2) e_3(\xi_3) \;,
\end{align*}
where we have used \eqref{eq:expan_diff} repeatedly to convert the derivatives of the nodal functions to edge functions.
Note that the divergence of the stress tensor and the body force, \eqref{eq:vol_expan}, are expanded in the same basis
and can therefore be added to represent (2a) in the finite-dimensional setting by
\begin{align}
\bm{D} \bm{\sigma}^{s,h}_m + f_m^{s,h}= \sum_{i,j,k=1}^N & \left [
\left ( T_{1m}^s \right )_{i,j,k} - \left ( T_{1m}^s \right )_{i-1,j,k} +
\left ( T_{2m}^s \right )_{i,j,k} - \left ( T_{2m}^s \right )_{i,j-1,k} + \right . \nonumber \\
& \left . \left ( T_{3m}^s \right )_{i,j,k} - \left ( T_{3m}^s \right )_{i,j,k-1} + \left (F_m^s \right ) \right ]
e_i(\xi_1) e_j(\xi_2) e_3(\xi_3) = 0 \;.
\label{eq:discretized_momentum_equation}
\end{align}
The important thing to note now, is that the $e_i(\xi_1) e_j(\xi_2) e_3(\xi_3)$ for $i,j,k=1,\ldots,N$ form
a {\em basis}, which means that
\[ \sum_{i,j,k=1}^N a_{i,j,k}e_i(\xi_1) e_j(\xi_2) e_3(\xi_3) = 0 \quad \Longleftrightarrow \quad
a_{i,j,k}=0 \;.\]
If we apply the linear independence of the basis functions to \eqref{eq:discretized_momentum_equation},
we see that we satisfy \eqref{eq:discretized_momentum_equation} point-wise, independent of
$(\xi_1,\xi_2,\xi_3)$, if and only if
\[ \left ( T_{1m}^s \right )_{i,j,k} - \left ( T_{1m}^s \right )_{i-1,j,k} +
\left ( T_{2m}^s \right )_{i,j,k} - \left ( T_{2m}^s \right )_{i,j-1,k} +
\left ( T_{3m}^s \right )_{i,j,k} - \left ( T_{3m}^s \right )_{i,j,k-1} + \left (F_m^s \right ) = 0\;,\]
which essentially is \eqref{eq:force_eq_C}.

In order to see that the expansion coefficients in the stress
representation, \eqref{eq:stress_expan}, actually represent the tractions over surfaces we can integrate, for instance,
$\sigma_{1m}^{s,h}$ over the a surface in the $(\xi_2,\xi_3)$-plane
\[ \int_{(\xi_2)_{j-1}}^{(\xi_2)_j} \int_{(\xi_3)_{k-1}}^{(\xi_3)_k} \left .
\sigma_{1m}^{s,h} \right |_{\xi_1=x_i} \,\mathrm{d}\xi_2 \mathrm{d} \xi_3 = \left (T_{1m}^s \right )_{i,j,k} \;,\]
where we used the Kronecker delta property of the Lagrange polynomials \eqref{eq:Kronecker_h} and the integral property of the edge polynomials \eqref{eq:Kronecker_e}.

The connection between elements is established by imposing that two neighboring elements share the degrees of freedom $T_{nm}^s$, i.e. the integrated traction or surface force at the interface. Although this seems to suggest that only the resulting surface forces are equilibrated, the surface forces are the expansion coefficients of polynomials of degree $N-1$ in both coordinates in the plane. Since we have $N$ surfaces in each coordinate direction, we completely fix the polynomials on both sides of the interface, which establishes a strong form of codiffusivity.

Note that the body force representation \eqref{eq:vol_expan} is an interpolation of the exact body force. We can exactly, point-wise, satisfy translational equilibrium if the exact body force has a piece-wise polynomial representation. This is, for instance, the case for the trivial body force $f=0$.

For the case where the body force is not a piecewise polynomial, (15) will be an approximation to the true body force and we satisfy point-wise translational equilibrium with respect to the approximate body force (15).
Note that if $\bm{f}$ is the exact body force density and $\bm{f}^h$ the polynomial approximation used in \eqref{eq:vol_expan}, then we still have that
\[ \int_{(\xi_1)_{i-1}}^{(\xi_1)_{i}} \int_{(\xi_2)_{j-1}}^{(\xi_2)_{j}} \int_{(\xi_3)_{k-1}}^{(\xi_3)_{k}} f_m = \int_{(\xi_1)_{i-1}}^{(\xi_1)_{i}} \int_{(\xi_2)_{j-1}}^{(\xi_2)_{j}} \int_{(\xi_3)_{k-1}}^{(\xi_3)_{k}} f_m^h \;,\]
for $i,j,k=1,\ldots,N$ and $m=1,2,3$ as a result of \eqref{eq:Kronecker_e}. So despite the fact that $\bm{f}^h$ is only a polynomial approximation of $\bm{f}$, the integrated $\bm{f}^h$ equals the integrated $\bm{f}$.

\section{Results I} \label{sec:Results_1}
A manufactured numerical test case is now presented. The test case is in 2D and Hooke's generalized law is the constitutive relation with the assumption of plane stress. This means that the compliance matrix is reduced to
\begin{equation*}
\bm{C} = \frac{1}{E} \begin{bmatrix} 1 & 0 & 0 & -\nu \\ 0 & (1+\nu) & 0 & 0 \\ 0 & 0 & (1+\nu) & 0 \\ -\nu & 0 & 0 & 1 \end{bmatrix} \;,
\end{equation*}
where $E$ is Young's modulus, which is set to $E=1$, and $\nu$ is Poisson's ratio, which is set to $\nu=0.3$.

\begin{figure}
\centering
\includegraphics[width=1.00\textwidth]{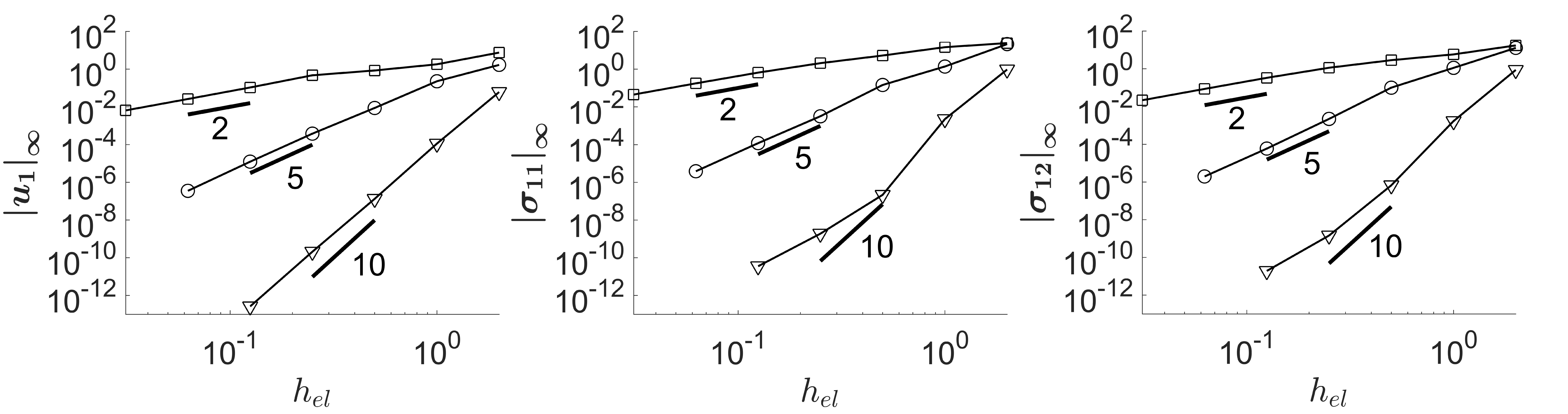} %
\caption{Convergence for $u_{1}$, $\sigma_{11}$ and $\sigma_{12}$ with respect to the element size $h_{el}$. $\Box:N=2$, $\circ:N=5$, $\bigtriangledown:N=10$.}
\label{fig:convergence_undef}
\end{figure}

\begin{figure}
\centering
\includegraphics[width=0.50\textwidth]{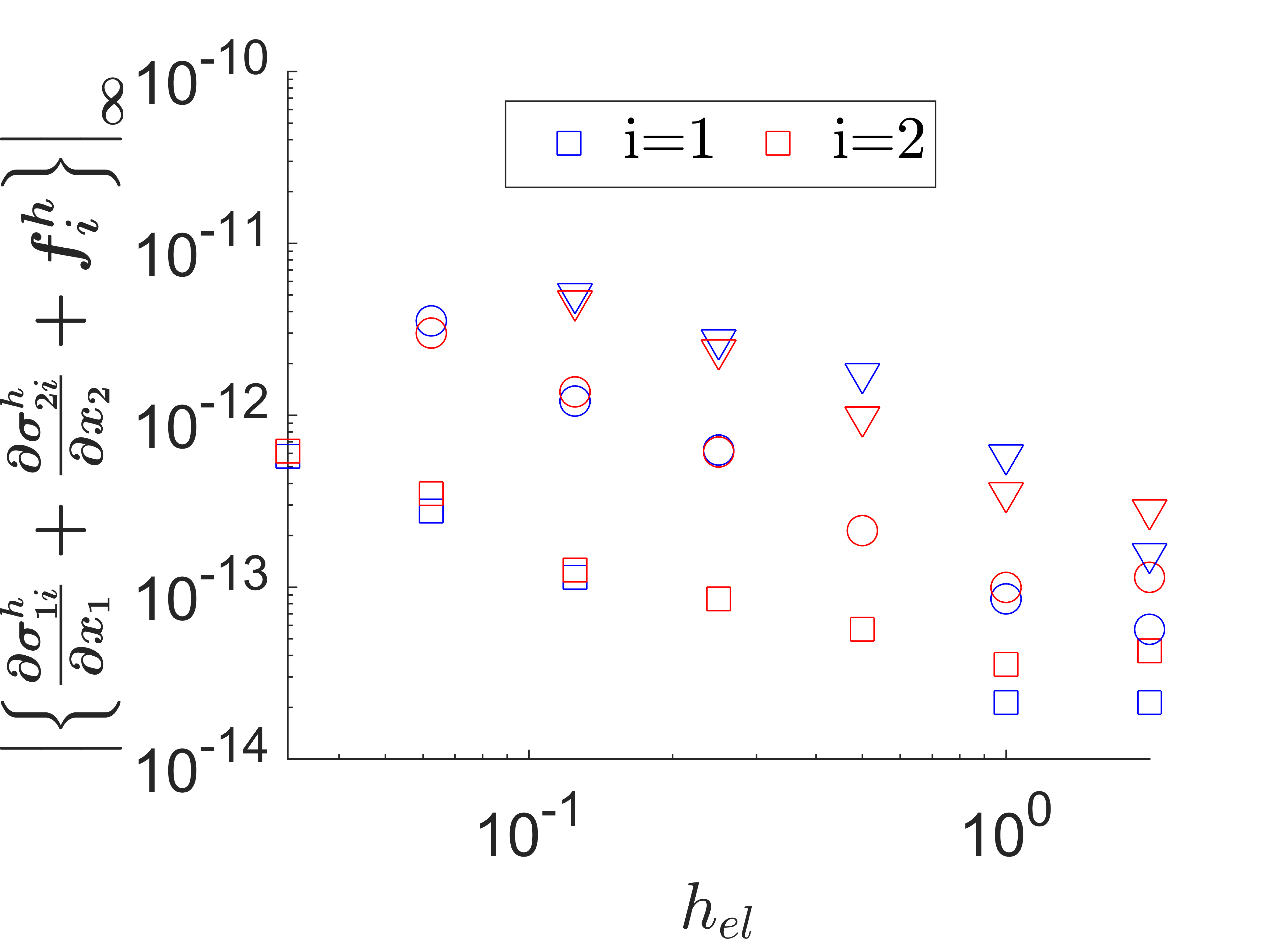}
\caption{The residual of the force equilibrium equations with respect to the element size $h_{el}$. $\Box:N=2$, $\circ:N=5$, $\bigtriangledown:N=10$.}
\label{fig:residual_force}
\end{figure}

The manufactured solution is given by
\begin{align*}
u_{1}(x_1,x_2) =& \sin(2 \pi x_1) \cos(2 \pi x_2) \;, \\
u_{2}(x_1,x_2) =& \cos(2 \pi x_1) \sin(2 \pi x_2) \;, \\
\sigma_{11}(x_1,x_2) = \sigma_{22}(x_1,x_2) =& \frac{\cos(2 \pi x_1) \cos(2 \pi x_2) 2 E \pi}{1-\nu} \;, \\
\sigma_{12}(x_1,x_2) =& -\frac{\sin(2 \pi x_1) \sin(2 \pi x_2) 2 E \pi}{1+\nu} \;, \\
f_{1}(x_1,x_2) =& -\frac{\sin(2 \pi x_1) \cos(2 \pi x_2) 8 E \pi^2}{1-\nu^2} \;, \\
f_{2}(x_1,x_2) =& -\frac{\cos(2 \pi x_1) \sin(2 \pi x_2) 8 E \pi^2}{1-\nu^2} \;,
\end{align*}
on the domain $\Omega \in [-1,1]^2$, and a displacement boundary condition is assigned to all the boundaries. A convergence study is  performed for the polynomials $N=\left\lbrace 2,5,10 \right\rbrace$. The convergence is evaluated in the infinity norm, i.e the maximum error evaluated in a preselected amount of points, where we have chosen $100 \times 100$ equispaced points in each element. The convergence is plotted in a logarithmic scale with respect to the size of the element, $h_{el}$. Convergence of the displacement field, normal stress field and shear stress field are presented in Figure~\ref{fig:convergence_undef}.
For a problem with smooth solution with polynomial degree $P$ we expect a convergence rate of $\mathcal{O}(h^{P+1})$ \cite{Karniadakis_2005}. Choosing $N$ to be the polynomial degree of the GLL Lagrange polynomials then it is seen that the highest polynomial order in \eqref{eq:disp_expan} and \eqref{eq:stress_expan} is $N-1$, and therefore the expected convergence rate for displacements and stresses are of $\mathcal{O}(h^N)$. So the expected slopes in the plots are $2$, $5$ and $10$, which are drawn as a reference. In general all the expanded variables show optimal convergence rates. More interesting is the plot of the residual of the force equilibrium equations in Figure~\ref{fig:residual_force}, which are calculated as
\begin{equation*}
R_{fe}(\bm{x}) = \left( \bm{\Psi}_3^G(\bm{x}) \bm{\mathcal{D}} \bm{\Delta}_{T} + \bm{\Psi}_3^G(\bm{x}) \bm{\Delta}_F  \right) \frac{1}{J(\bm{x})} \;.
\end{equation*}
The plot shows that the force equilibrium equations are satisfied to machine precision in a finite volume setting independent of the resolution of the computational domain. This means that the numerical errors are confined to the constitutive equations and the equilibrium of moments.

\vskip 0.3cm
\noindent
\textbf{Remark 4} In Figure~\ref{fig:residual_force} the $L^{\infty}$-error in conservation of translational equilibrium is satisfied up $O(10^{-14}-10^{-11})$, which is in the order of machine precision. If we refine the mesh, i.e. decrease $h$, or take a higher polynomial degree, increasing $N$, the ``error" in Figure~\ref{fig:residual_force} increases, which is counter-intuitive. The reason is that in both cases the number of degrees of freedom increases and as a result the condition number of the linear system grows leading to higher ``errors". If exact arithmetic were used, the error would be exactly zero and remain zero for all $h$ and $N$. Note that a similar error growth can be seen in Figure~\ref{fig:residual_force_def}.

\section{Transformations} \label{sec:transformations}
\begin{figure}
\centering
\includegraphics[width=1.0\textwidth]{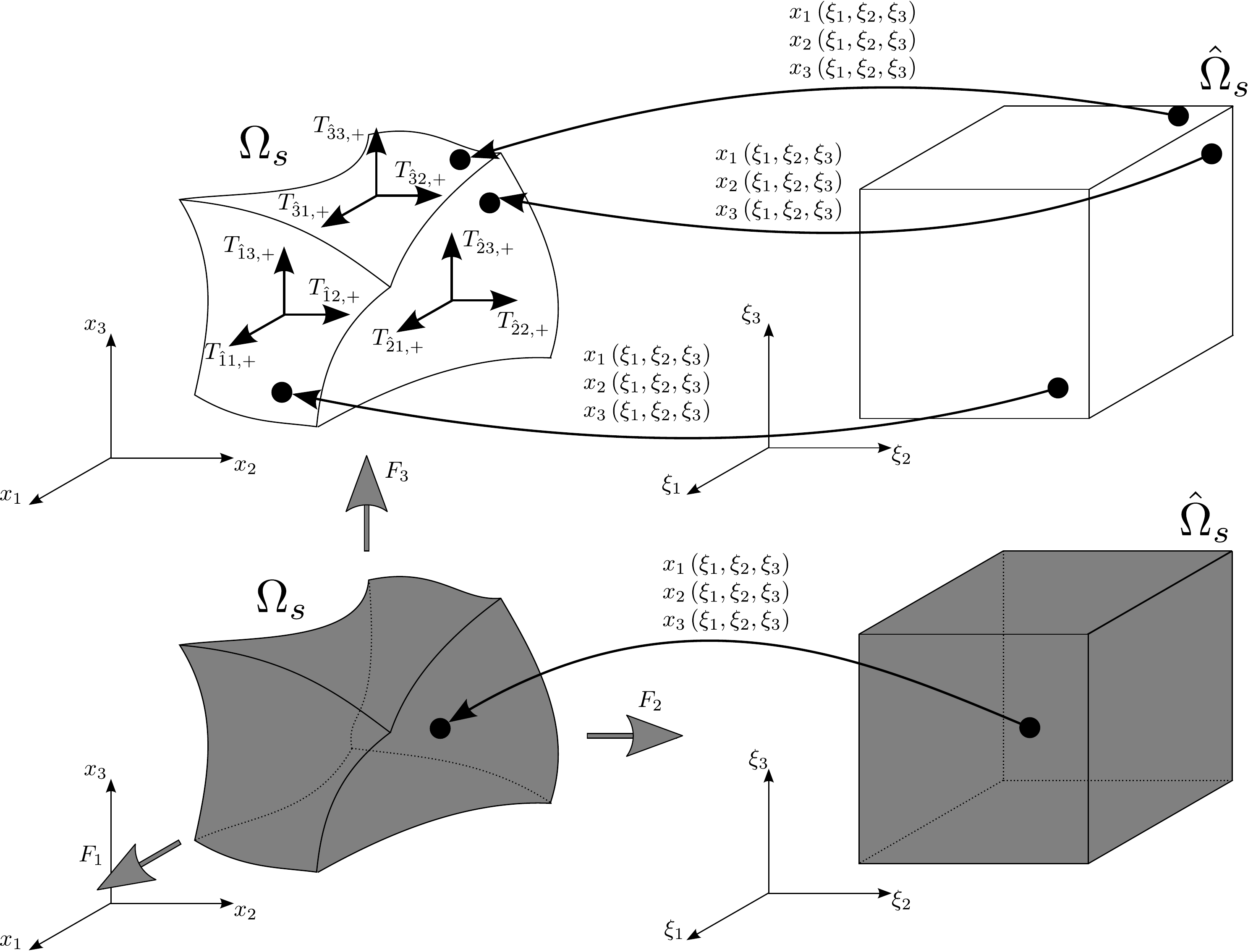}
\caption{The forces on a deformed element, $\Omega_s$, in $\mathbb{R}^3$, and the mapping from a reference element, $\hat{\Omega}_s$.}
\label{fig:Forces}
\end{figure}

For practical applications we need to expand the aforementioned approach to deformed grids, and we let $\Omega$ be a deformed domain and $\Omega_s$, depicted in Figure~\ref{fig:Forces}, be a sub-domain or element of $\Omega$, and $\hat{\Omega}_s$ is a reference element with a mapping to $\Omega_s$.
The geometrical association of the surface force components and body force components can be expressed with respect to the reference element, but the direction of the components is maintained in the physical reference frame. The surface force components are now denoted by $T_{\hat{\i} j}$, where the index $\hat{\i}$ is the reference Cartesian direction of the boundary surface of $\hat{\Omega}_s$, while the $j$ index denotes the direction of the surface force component with respect to the $x_j$ basis. Arranging the individual surface force components $T_{\hat{\i}j}$ in the column vector $\bm{\Delta}_{\hat{T}}$, then the equilibrium of forces is given by
\begin{equation}
\bm{\mathcal{D}} \bm{\Delta}_{\hat{T}} = - \bm{\Delta}_{\hat{F}} \;, \label{eq:force_ref}
\end{equation}
where $\bm{\mathcal{D}}$ is the same matrix as in \eqref{eq:force_phys}, and $\bm{\Delta}_{\hat{F}}$ is a column vector with all components of the body forces of the deformed elements. The connection between stress components and surface force components is now given by
\begin{align*}
T_{\hat{1}j} =& \iint \sigma_{\hat{1}j} \, \mathrm{d}\xi_2 \mathrm{d}\xi_3 \\
T_{\hat{2}j} =& \iint \sigma_{\hat{2}j} \, \mathrm{d}\xi_3 \mathrm{d}\xi_1 \\
T_{\hat{3}j} =& \iint \sigma_{\hat{3}j} \, \mathrm{d}\xi_1 \mathrm{d}\xi_2 \;,
\end{align*}
while the body force components are given by
\begin{equation*}
F_{j} = \iiint_{\hat{\Omega}_s} \hat{f}_j \, \mathrm{d}\xi_1 \mathrm{d}\xi_2 \mathrm{d}\xi_3 \;,
\end{equation*}
where $\hat{f}_j$ are the body force components associated to the reference element. Note that $\sigma_{\hat{\i}j}$ resembles a component of the first Piola--Kirchhoff stress tensor, \cite[p.164]{Reddy_2013}; however, where the first Piola--Kirchhoff stress tensor is defined in the undeformed geometry, the mapping employed here uses the geometry of the reference element. A similar approach was used in \cite{Santos_2014}, where hybrid stress elements were derived based on transformations between the Cauchy stress tensor and stress tensors resembling the first and second Piola--Kirchhoff stress tensors. As seen in \eqref{eq:force_ref} the discrete version of the force equilibrium equations is free of interpolations also in the case of deformed elements. We saw in Section~\ref{sec:discussion} that we can satisfy \eqref{eq:discretized_momentum_equation} in the strong sense, when we satisfy \eqref{eq:force_ref}. Any transformation that maps the basis functions $e_i(\xi_1)e_j(\xi_2)e_k(\xi_3)$ to basis functions which are still linear independent \MIG{implies} that the expansion coefficients need to vanish in order to satisfy the equation. The expansion coefficients remain unchanged under the mapping, therefore we still need \eqref{eq:force_ref} to hold. The key issue is that {\em all terms in \eqref{eq:discretized_momentum_equation} are expanded in the same basis and therefore can be added. This property is a direct consequence of the degrees of freedom (Section~\ref{sec:eq_force}) and the basis functions (Section~\ref{sec:expan})}.

In \eqref{eq:stress_strain_rel_strain_exp} the constitutive equations are written with respect to the Cauchy stress components, and hence a mapping between $\sigma_{ij}$ and $\sigma_{\hat{m}n}$ must be introduced. According to \cite[p.165]{Reddy_2013} the relation between the Cauchy stress tensor, $\bm{\sigma}$, and the first Piola--Kirchhoff stress tensor, $\hat{\bm{\sigma}}$ is given by
\begin{equation*}
\sigma_{ij} = \frac{1}{J} \mathcal{F}_{ik} \hat{\sigma}_{kj} \;,
\end{equation*}
where $\mathcal{F}_{ik}$ are components of the deformation gradient matrix given by
\begin{equation*}
\bm{\mathcal{F}} = \begin{bmatrix}
\frac{\partial x_1}{\partial \xi_1} & \frac{\partial x_1}{\partial \xi_2} & \frac{\partial x_1}{\partial \xi_3} \\
\frac{\partial x_2}{\partial \xi_1} & \frac{\partial x_2}{\partial \xi_2} & \frac{\partial x_2}{\partial \xi_3} \\
\frac{\partial x_3}{\partial \xi_1} & \frac{\partial x_3}{\partial \xi_2} & \frac{\partial x_3}{\partial \xi_3}
\end{bmatrix} \;,
\end{equation*}
and $J = \det(\bm{\mathcal{F}})$ is the Jacobian. Note that normally $\bm{\mathcal{F}}$ is the connection between the undeformed and deformed material state, while in the current case it is a connection between the \MIG{actual geometry} of the element and the reference element. The relation can be formulated in engineering notation as
\begin{equation}
\bm{\sigma} = \frac{1}{J} \begin{bmatrix}
\bm{\mathcal{F}} & 0 & 0 \\
0 & \bm{\mathcal{F}} & 0 \\
0 & 0 & \bm{\mathcal{F}}
\end{bmatrix} \hat{\bm{\sigma}} = \frac{1}{J} \bm{\mathcal{F}}^e \hat{\bm{\sigma}} \;, \label{eq:Cauchy_f_Piola_rel}
\end{equation}
with
\begin{equation*}
\hat{\bm{\sigma}} = \begin{Bmatrix}
\sigma_{\hat{1}1} & \sigma_{\hat{2}1} & \sigma_{\hat{3}1} & \sigma_{\hat{2}1} & \sigma_{\hat{2}2} & \sigma_{\hat{2}3} & \sigma_{\hat{3}1} & \sigma_{\hat{3}2} & \sigma_{\hat{3}3}
\end{Bmatrix}^T \;.
\end{equation*}

A volume integration of an element in the physical reference frame can be performed through the reference element by
\begin{equation}
\iiint_{\Omega_s} \, \mathrm{d}x_1 \mathrm{d}x_2 \mathrm{d}x_3 = \iiint_{\hat{\Omega}_s} J \, \mathrm{d}\xi_1 \mathrm{d}\xi_2 \mathrm{d}\xi_3 \;, \label{eq:vol_int}
\end{equation}
and hence
\begin{equation*}
\hat{f}_j = f_j J \;.
\end{equation*}

Because we consider rotation as a nodally defined quantity, the expansion of the rotation components in \eqref{eq:expan_rot} are with respect to the element basis and a connection to the physical coordinate system is given by
\begin{equation}
\bm{\omega} = \hat{\bm{\omega}} \;. \label{eq:rot_trans}
\end{equation}

Equation \eqref{eq:Cauchy_f_Piola_rel}, \eqref{eq:vol_int} and \eqref{eq:rot_trans} are inserted into \eqref{eq:approx_eq_sum} and yield
\begin{equation}
\begin{bmatrix}
\bm{H}^h & \left ( \bm{V}^h \bm{\mathcal{D}} \right )^T & -\left( \bm{R}^h \right)^T  \\
\bm{V}^h \bm{\mathcal{D}}  & \bm{0} & \bm{0} \\
-\bm{R}^h & \bm{0} &  \bm{0}
\end{bmatrix} \begin{Bmatrix}
\bm{\Delta}_{\hat{T}} \\  \bm{\Delta}_u \\ \bm{\Delta}_{\hat{O}}
\end{Bmatrix} = \begin{Bmatrix}
\bm{B}^h \bm{\Delta}_{\bar{u}} \\ - \bm{V}^h \bm{\Delta}_{\hat{F}} \\
\bm{0}
\end{Bmatrix} \;, \label{eq:lin_eq_sys}
\end{equation}
is produced with
\begin{subequations} \label{eq:lin_eq_sys_sub_mat}
\begin{align}
\bm{V}^h =& \sum\limits_{s=1}^{N_e} \sum\limits_{i=0}^{N-1} \sum\limits_{j=0}^{N-1} \sum\limits_{k=0}^{N-1} \left( \left( \bm{\Psi}_0^s(\tilde{\bm{\xi}}_{i,j,k}) \right)^T \bm{\Psi}_3^s(\tilde{\bm{\xi}}_{i,j,k}) \tilde{w}_i \tilde{w}_j \tilde{w}_k \right) \;, \\
\bm{H}^h =& \frac{1}{J} \sum\limits_{s=1}^{N_e} \sum\limits_{i=0}^{N} \sum\limits_{j=0}^{N} \sum\limits_{k=0}^{N} \left( \left( \bm{\Psi}_2^s(\bm{\xi}_{i,j,k}) \right)^T \left( \bm{\mathcal{F}}^e \right)^T \bm{C} \bm{\mathcal{F}}^e \bm{\Psi}_2^s(\bm{\xi}_{i,j,k}) w_i w_j w_k \right) \;, \\
\bm{R}^h =& \sum\limits_{s=1}^{N_e} \sum\limits_{i=0}^{N} \sum\limits_{j=0}^{N} \sum\limits_{k=0}^{N} \left( \left( \bm{\Psi}_{0}^s(\bm{\xi}_{i,j,k}) \right)^T
\bm{R} \bm{\mathcal{F}}^e \tilde{\bm{\Psi}}_2^s(\bm{\xi}_{i,j,k}) w_i w_j w_k \right) \;, \\
\bm{B}^h =& \sum\limits_{s=1}^{N_{BCu}} \sum\limits_{i=0}^{N} \sum\limits_{j=0}^{N} \left( \left( \bm{\Psi}_2^s(\tilde{\bm{\xi}}_{i,j,k}) \right)^T \tilde{\bm{\Psi}}_0^s(\tilde{\bm{\xi}}_{i,j,k}) \tilde{w}_i \tilde{w}_j \right) \;.
\end{align}
\end{subequations}

\MIG{The linear system (\ref{eq:lin_eq_sys}) corresponds to a {\em mixed formulation}, \cite{BrezziFortin,BoffiBrezziFortin}. Well-posedness of such a system is ensured when $\bm{\mathcal{D}}$ and $\bm{R}^h$ are both surjective maps and when the symmetric sub-matrix $\bm{H}^h$ is coercive on the null space of 
\[ \left [ \begin{array}{c}
\bm{V}^h \bm{\mathcal{D}} \\
-\bm{R}^h
\end{array} \right ]\;,\]
By construction, see Section \ref{sec:eq_force}, the operator $\bm{V}^h \bm{\mathcal{D}}$ is surjective. The operator $\bm{R}^h$ is surjective if we define the degrees of freedom for the rotation $\omega$ in the Gauss points, but not when rotation is defined in the Gauss-Lobatto points as will be shown in Section~\ref{sec:Comp_energy}. The matrix $\bm{H}^h$ is $L^2$-coercive and for those stress fields that satisfy $\bm{\mathcal{D}} \Delta_{\hat{T}}=0$ the discrete operator $\bm{H}^h$ is even coercive on the space
\[ H_\sigma(\bm{\mathcal{D}}):= \left \{ \sigma_{ij} \in L^2(\Omega) \,\big | \, \bm{\mathcal{D}} \sigma \in \left [ L^2(\Omega) \right ]^2 \, \right \} \;. \]}

\section{Results II} \label{sec:Results_2}
Some numerical results involving non-orthogonal domains are now presented. First the manufactured solution from Section~\ref{sec:Results_1} is extended to distorted elements, and then the stress concentration around a circular hole is investigated.

\subsection{Manufactured solution}
\begin{figure}
\centering
\includegraphics[width=1.00\textwidth]{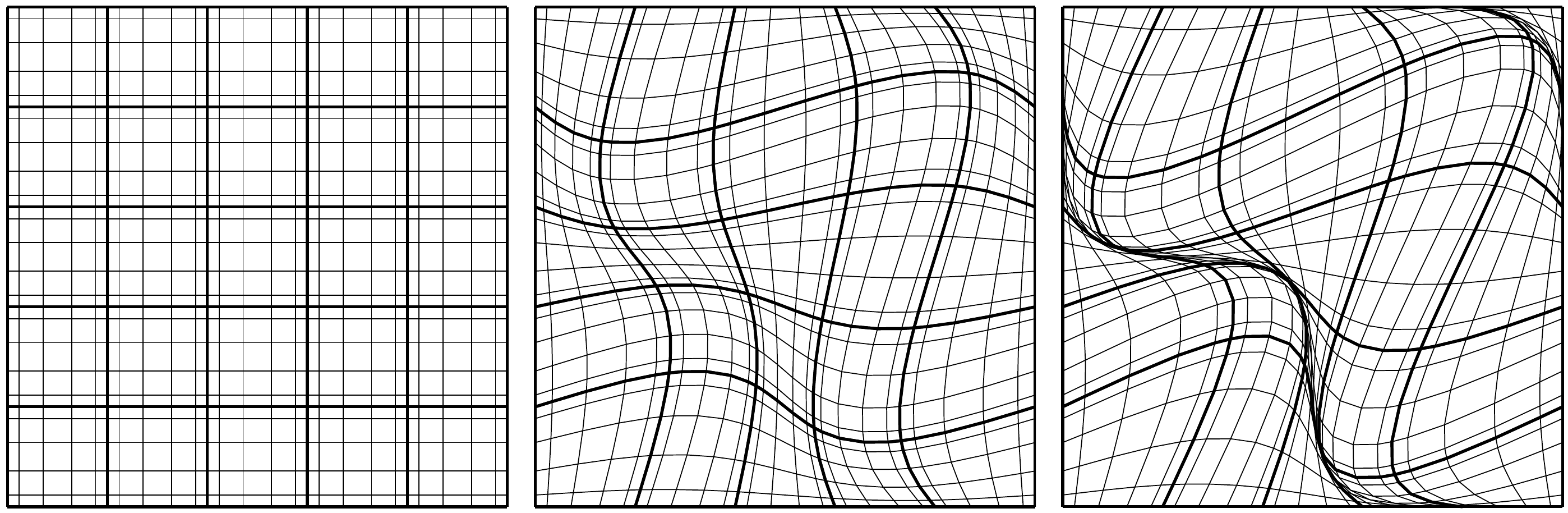}
\caption{Deformed grids for $5 \times 5$ elements with $N=5$ and $c = \left\lbrace 0,0.15,0.3 \right\rbrace$.}
\label{fig:mesh_3}
\end{figure}

The mapping \cite{Kreeft_stokes,Gerritsma_2012}
\begin{subequations} \label{eq:mapping}
\begin{align}
x_1(\xi_1,\xi_2) =& \xi_1 + c \sin(\pi \xi_1) \sin(\pi \xi_2) \;, \\
x_2(\xi_1,\xi_2) =& \xi_2 + c \sin(\pi \xi_1) \sin(\pi \xi_2) \;,
\end{align}
\end{subequations}
is now introduced to the manufactured solution from Section~\ref{sec:Results_1}, which is shown in Figure~\ref{fig:mesh_3} for $c = \left\lbrace 0,0.15,0.3 \right\rbrace$.
Note that the grid for $c=0$ coincides with the grid used in Section~\ref{sec:Results_1}.

\begin{figure}
\centering
\includegraphics[width=1.00\textwidth]{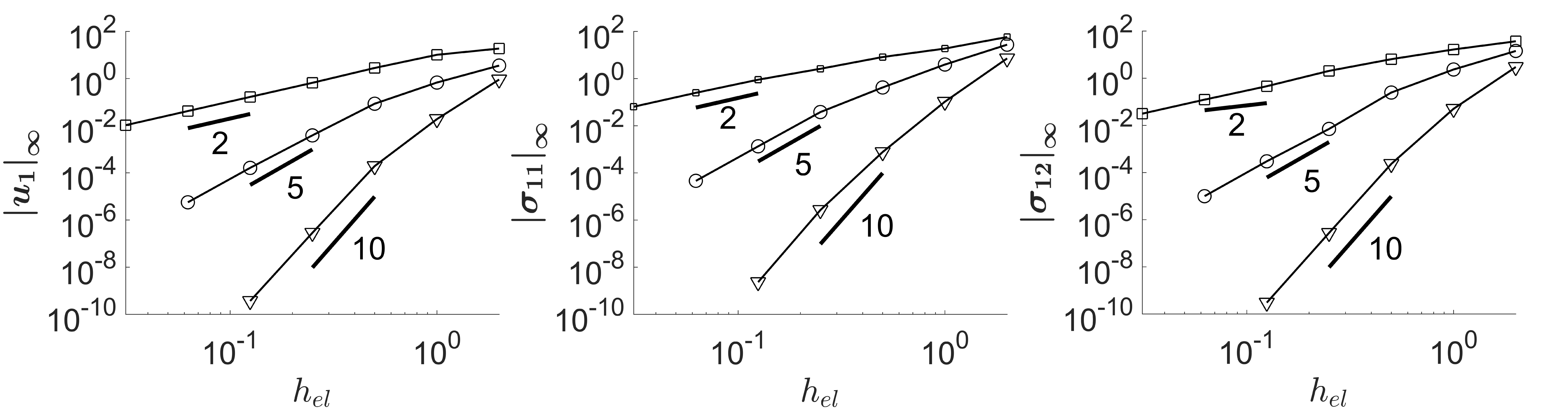} %
\caption{Convergence for $u_{1}$, $\sigma_{11}$ and $\sigma_{12}$ with respect to the undeformed element size $h_{el}$ using the mapping in \eqref{eq:mapping} with $c=0.15$. $\Box:N=2$, $\circ:N=5$, $\bigtriangledown:N=10$.}
\label{fig:convergence_c0_15}
\end{figure}

\begin{figure}
\centering
\includegraphics[width=1.00\textwidth]{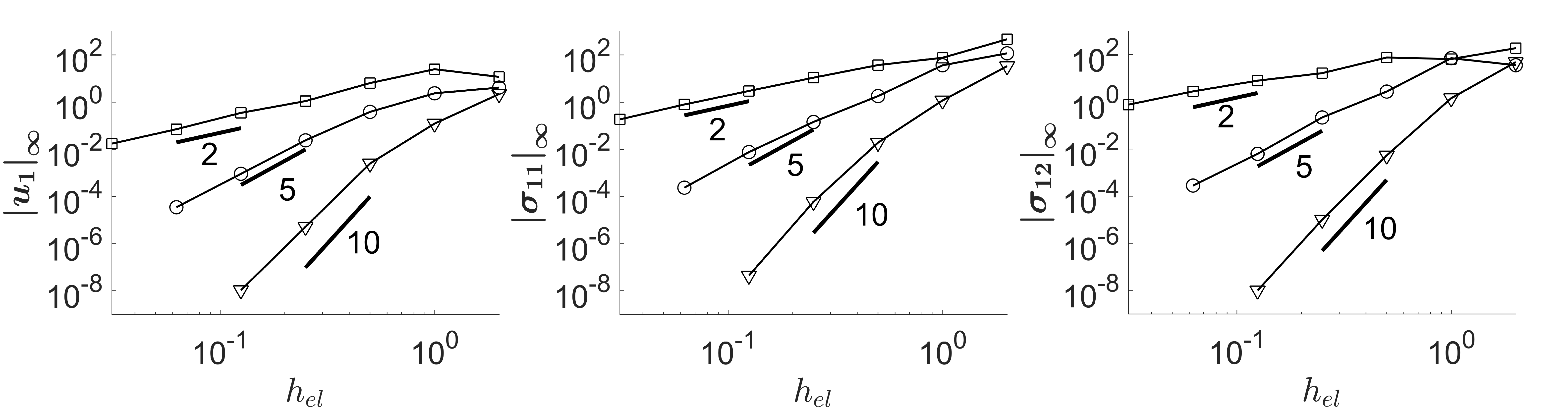} %
\caption{Convergence for $u_{1}$, $\sigma_{11}$ and $\sigma_{12}$ with respect to the undeformed element size $h_{el}$ using the mapping in \eqref{eq:mapping} with $c=0.3$. $\Box:N=2$, $\circ:N=5$, $\bigtriangledown:N=10$.}
\label{fig:convergence_c0_30}
\end{figure}

\begin{figure}
\centering
\includegraphics[width=1.00\textwidth]{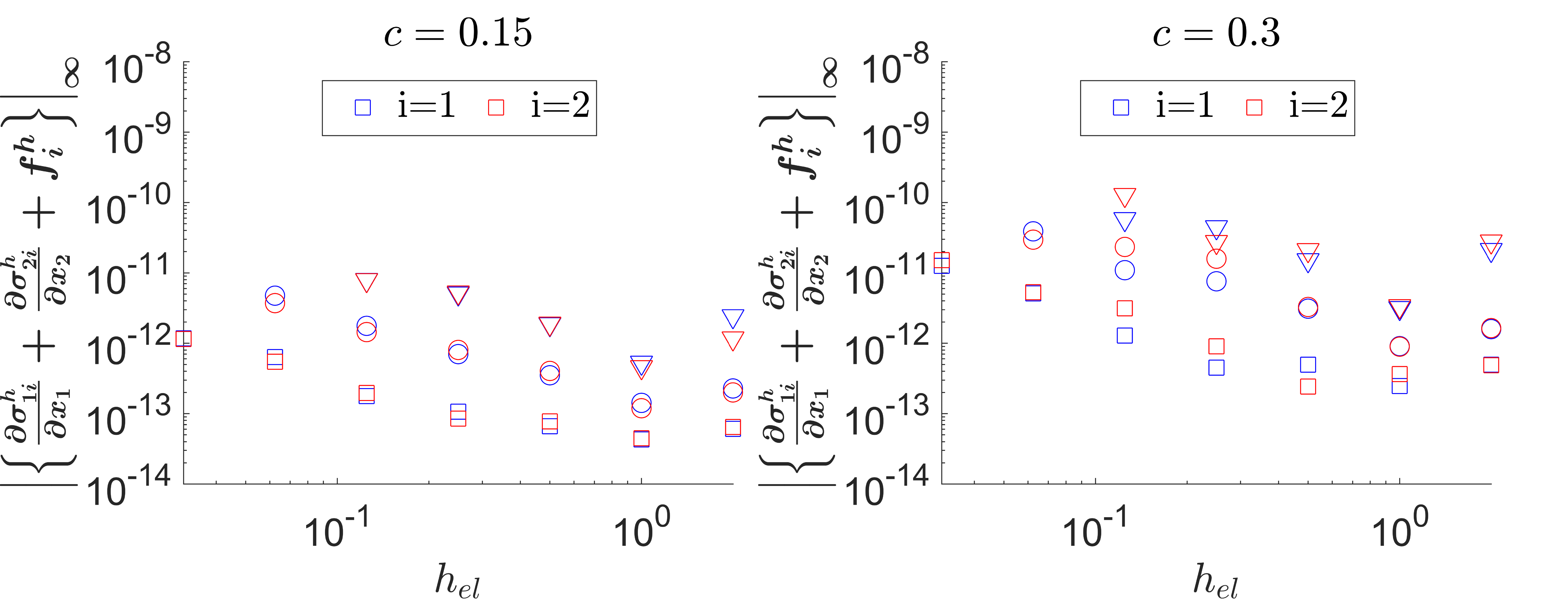}
\caption{The residual of the force equilibrium equations with respect to the undeformed element size $h_{el}$ for $c = \left\lbrace 0.15,0.3 \right\rbrace$. $\Box:N=2$, $\circ:N=5$, $\bigtriangledown:N=10$.}
\label{fig:residual_force_def}
\end{figure}

The same convergence study as in Section~\ref{sec:Results_1} is performed, but with $c=0.15$ and $c=0.3$. As seen in Figure~\ref{fig:convergence_c0_15} and Figure~\ref{fig:convergence_c0_30} we get the expected convergence rates.
The force equilibrium equations are also satisfied to machine precision in a finite volume setting for both $c=0.15$ and $c=0.3$, as shown in Figure~\ref{fig:residual_force_def}.
On curvilinear grids the integrands are no longer polynomial and therefore Gauss or Gauss--Lobatto integration is no longer exact. In all calculations the integrals were also evaluated with very high order Gauss/Gauss--Lobatto integration. The integration error is much smaller than the numerical error and high order integration does not change the results appreciably.

\subsection{Plate with Circular Hole}
\begin{figure}
\centering
\includegraphics[width=1.00\textwidth]{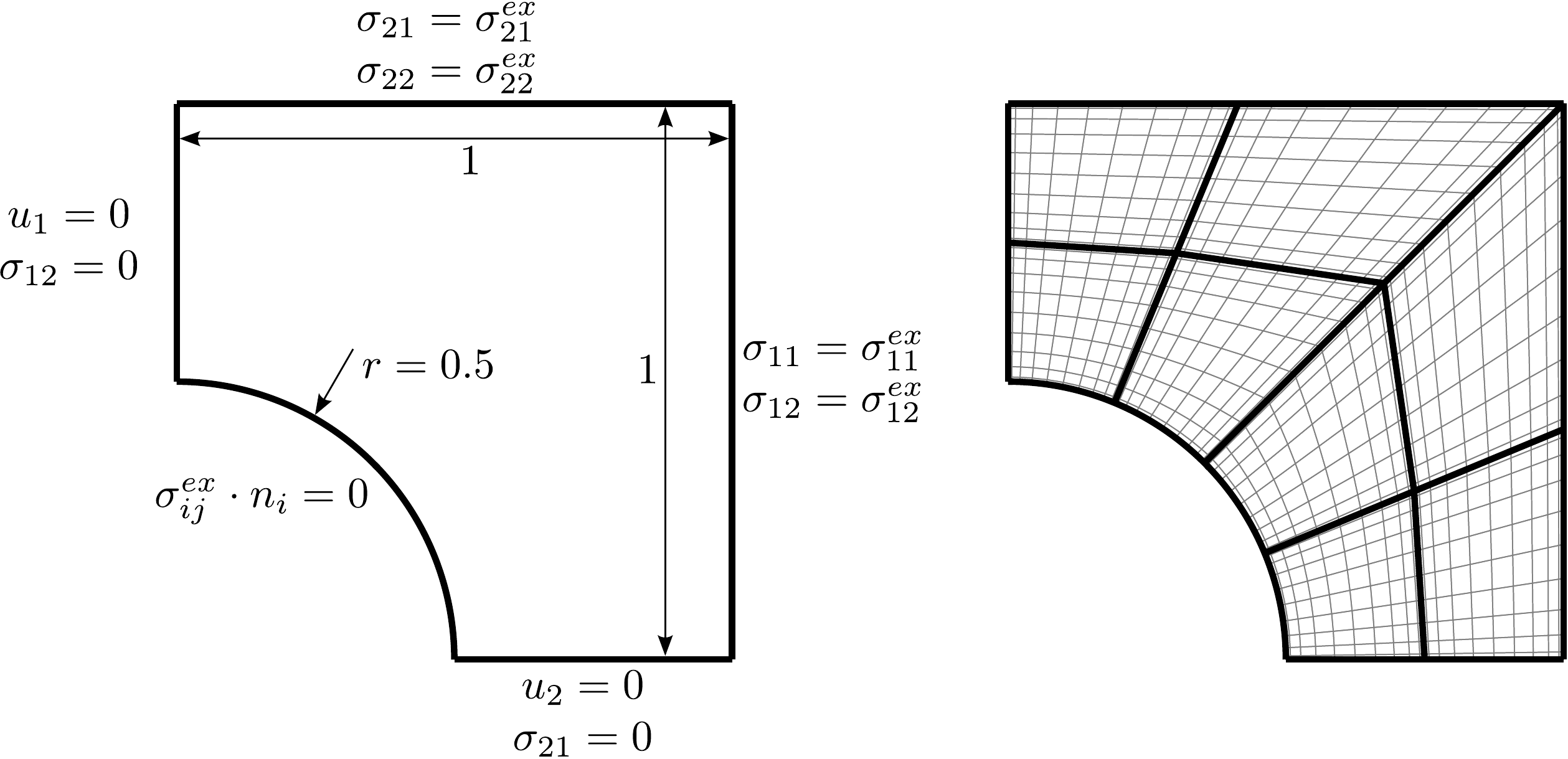} %
\caption{The quarter section of a circular hole of size  $r=0.5$ in a plate of size $2 \times 2$. Left: The boundary conditions. Right: The grid consisting of 8 spectral elements with polynomials degree $N=10$.}
\label{fig:circ_hole_test_case}
\end{figure}

\begin{figure}
\centering
\includegraphics[width=0.50\textwidth]{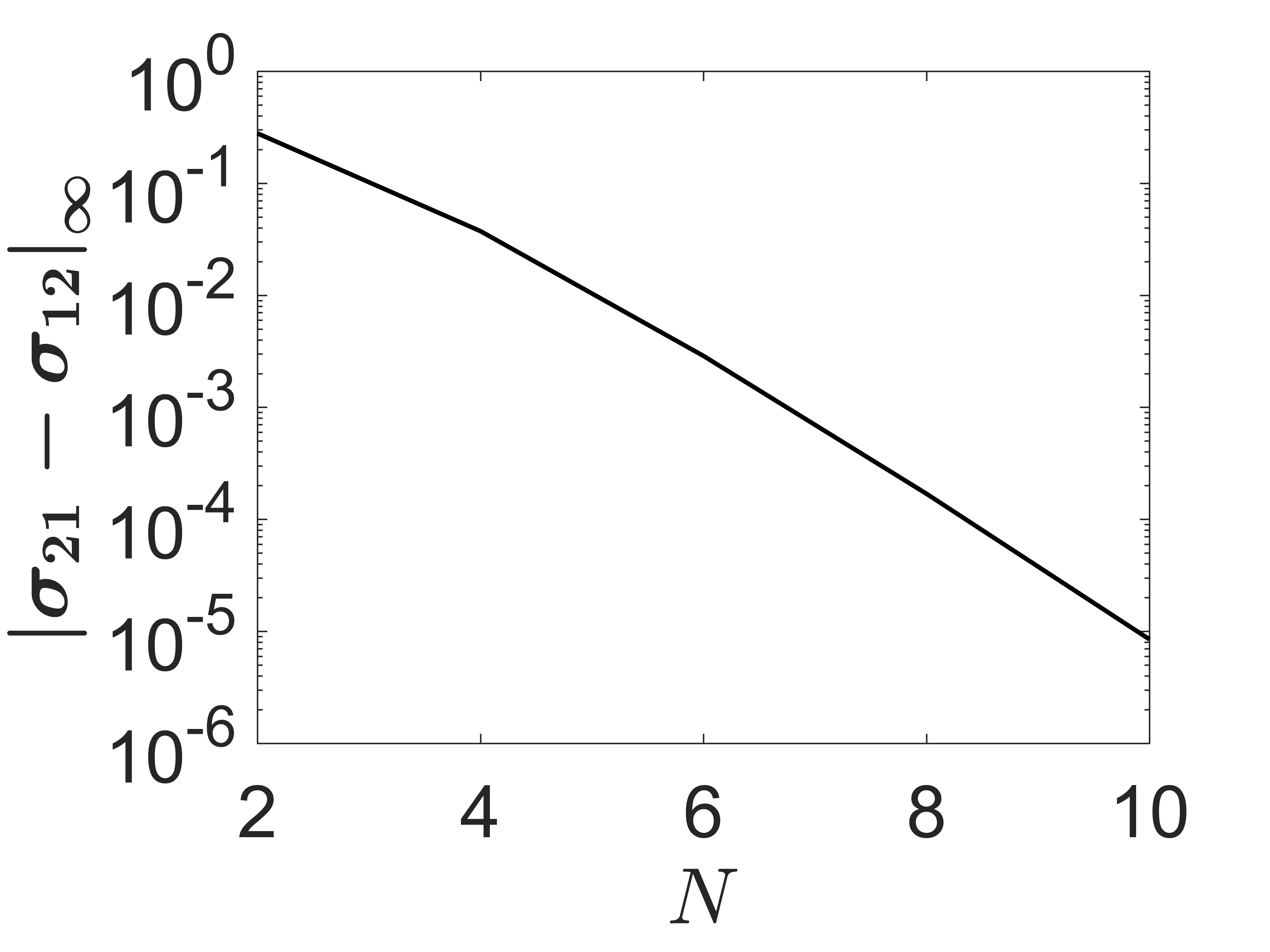}
\caption{The difference between the shear stress components using the grid in Figure~\ref{fig:circ_hole_test_case} with $N=\left\lbrace 2,4,6,8,10 \right\rbrace$.}
\label{fig:circ_dts}
\end{figure}

\begin{figure}
\centering
\includegraphics[width=1.0\textwidth]{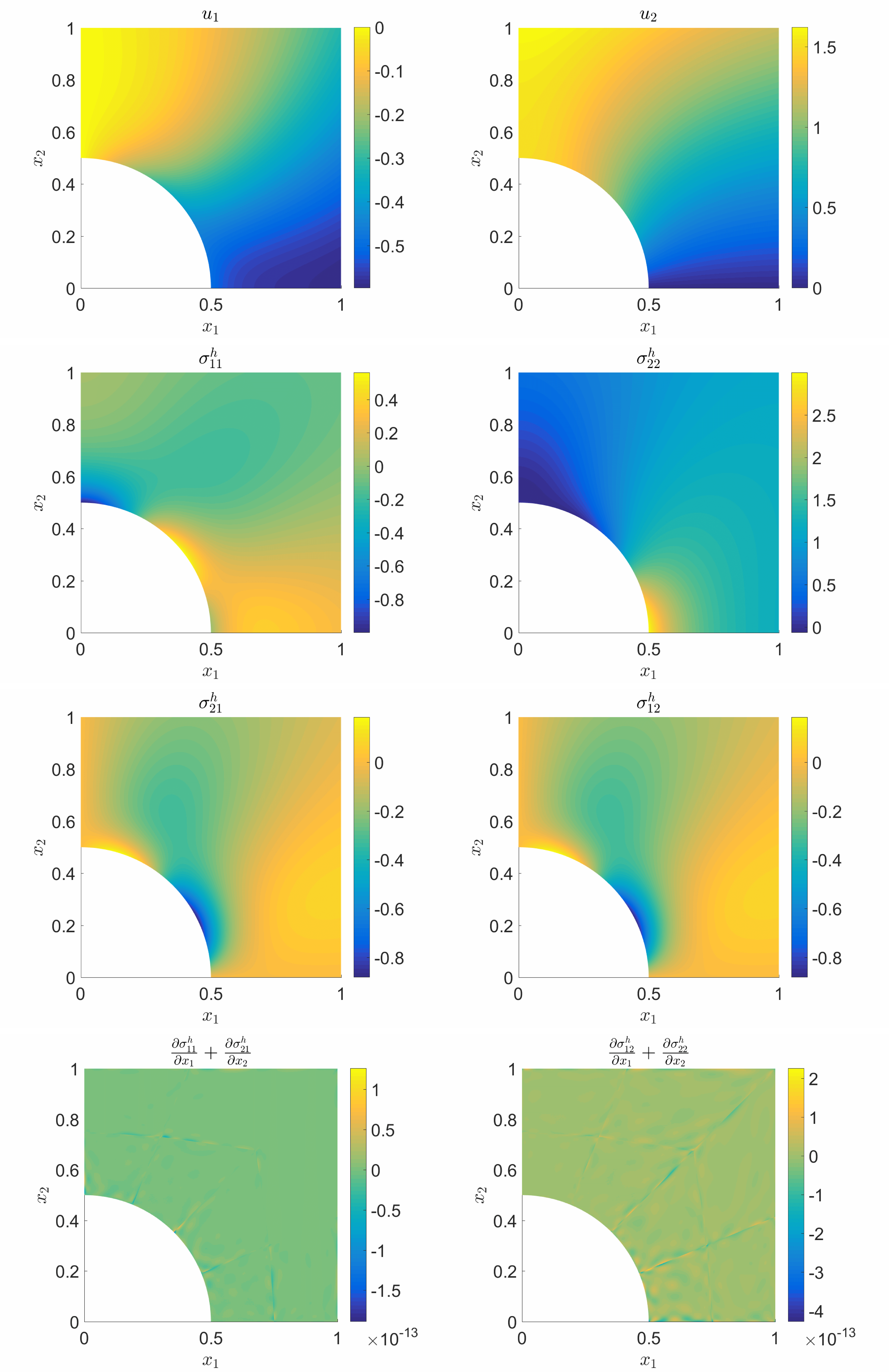}
\caption{The solution of the test case in Figure~\ref{fig:circ_hole_test_case} as well as the residual of the force equilibrium equations.}
\label{fig:circ_hole_test_case_sol}
\end{figure}

This test case is for the stress concentration of a circular hole in the plate shown in Figure~\ref{fig:circ_hole_test_case}, which has a unidirectional loading. The exact solution for a circular hole in an infinite large plate is given in many textbooks on stress analysis (e.g., see \cite[p. 306]{Reddy_2013}). By using symmetry conditions only a quarter of the domain is considered. The boundary condition and the grid resolution are shown in Figure~\ref{fig:circ_hole_test_case}.
Note that the traction boundary conditions (BCs) are implemented strongly through surface forces on $\Gamma_t$. The elements are mapped to the reference domain with a transfinite map, \cite{Gordon_1973}. The solution of the displacement field, stress field  and the residual of the force equilibrium equations are plotted in Figure~\ref{fig:circ_hole_test_case_sol}, and the residual of the force equilibrium equations are of order $10^{-13}$. Note that since we have no body forces the equilibrium equations are not only satisfied in a finite volume setting but they are satisfied in a general sense. The largest errors in the domain are listed in Table~\ref{tab:circ_res}. As seen from the norms the solution of $\sigma_{1 2}^h$ and $\sigma_{2 1}^h$ are not equal, which is expected as they are expanded by different polynomial see \eqref{eq:stress_expan}. To see how the difference between the two shear stress components decreases as the resolution increases, a convergence study is performed. The grid with $8$ spectral elements in Figure~\ref{fig:circ_hole_test_case} is maintained, and $N$ is varied from $N=2$ to $N=10$. The result is plotted in Figure~\ref{fig:circ_dts}.
It is observed that the difference decreases with an exponential rate.

\begin{table}
\caption{The infinite norm for the error of the displacement and stress fields in the test case of the circular hole in Figure~\ref{fig:circ_hole_test_case}, and $(x_{1,i},x_{2,i})$ are the coordinates of an amount of preselected points.}
\centering
\tabsize
\begin{tabular}{cc}
\toprule
$z_i$ & $\left| \bm{z} \right|_\infty$ \\
\midrule
$u_1^h(x_{1,i},x_{2,i}) - u_1^{ex}(x_{1,i},x_{2,i})$ & $5.4547 \cdot 10^{-7}$ \\
$u_2^h(x_{1,i},x_{2,i}) - u_2^{ex}(x_{1,i},x_{2,i})$ & $5.7689 \cdot 10^{-7}$ \\
$\sigma_{11}^h(x_{1,i},x_{2,i}) - \sigma_{11}^{ex}(x_{1,i},x_{2,i})$ & $6.7320 \cdot 10^{-6}$ \\
$\sigma_{21}^h(x_{1,i},x_{2,i}) - \sigma_{21}^{ex}(x_{1,i},x_{2,i})$ & $6.0327 \cdot 10^{-6}$ \\
$\sigma_{12}^h(x_{1,i},x_{2,i}) - \sigma_{12}^{ex}(x_{1,i},x_{2,i})$ & $5.8757 \cdot 10^{-6}$ \\
$\sigma_{22}^h(x_{1,i},x_{2,i}) - \sigma_{22}^{ex}(x_{1,i},x_{2,i})$ & $6.6669 \cdot 10^{-6}$ \\
\bottomrule
\end{tabular}
\label{tab:circ_res}
\end{table}

\section{Energy considerations and equilibrium of moments} \label{sec:Comp_energy}
The complementary strain energy is given by, \cite{Veubeke_1964,Veubeke_1965}
\begin{equation}
U_C^h = \frac{1}{2} \int_{\Omega} \left( \bm{\sigma}^h \right)^T \bm{C} \bm{\sigma}^h \, \mathrm{d} \Omega \;. \label{eq:comp_strain}
\end{equation}
This energy will converge towards the exact solution monotonically from above in equilibrium models with homogeneous displacement boundary conditions. Following \cite{Almeida_Continuity,Veubeke_1964,Veubeke_1965,Debongnie_1995} the stress field is divided into a homogeneous stress field and a particular stress field, $\bm{\sigma}^p$. The particular stress field satisfies \eqref{eq:cons_lin_mom} a priori and the homogeneous part is divergence free. Hence the body force part of \eqref{eq:cons_lin_mom} is moved to \eqref{eq:stress_strain_rel_strain_exp}. We can implement this by changing the linear equation system in \eqref{eq:lin_eq_sys} to
\begin{equation*}
\begin{bmatrix}
\bm{H}^h & \left ( \bm{V}^h \bm{\mathcal{D}} \right )^T & -\left( \bm{R}^h \right)^T \\
\bm{V}^h \bm{\mathcal{D}}  & \bm{0} & \bm{0} \\
-\bm{R}^h & \bm{0} &  \bm{0}
\end{bmatrix} \begin{Bmatrix}
\bm{\Delta}_{\hat{T}} \\  \bm{\Delta}_u \\ \bm{\Delta}_{\hat{O}}
\end{Bmatrix} = \begin{Bmatrix}
\bm{B}^h \bm{\Delta}_{\bar{u}} - \bm{H}^h_p \bm{\Delta}_{\sigma_p} \\
\bm{0} \\
\bm{0}
\end{Bmatrix} \;,
\end{equation*}
where the sub-matrices are found in \eqref{eq:lin_eq_sys_sub_mat}, $\bm{H}^h_p$ is given by
\begin{equation*}
\bm{H}^h_p = \sum\limits_{s=1}^{N_e} \sum\limits_{i=0}^{N} \sum\limits_{j=0}^{N} \sum\limits_{k=0}^{N} \left( \left( \bm{\Psi}_2^s(\bm{\xi}_{i,j,k}) \right)^T \left( \bm{\mathcal{F}}^e \right)^T \bm{C} w_i w_j w_k \right) \;,
\end{equation*}
and $\bm{\Delta}_{\sigma_p}$ is a column vector with the particular stress field evaluated at the GLL point in all elements.

A manufactured solution with homogeneous displacement BCs is given by
\begin{subequations} \label{eq:comp_energy_man_sol}
\begin{align}
u_{1}(x_1,x_2) =& \sin(2 \pi x_1) \sin(2 \pi x_2) \;, \\
u_{2}(x_1,x_2) =& \sin(2 \pi x_1) \sin(2 \pi x_2) \;, \\
\sigma_{11}(x_1,x_2) =& 2 E \pi \frac{\cos(2 \pi x_1) \sin(2 \pi x_2) + \nu \sin(2 \pi x_1) \cos(2 \pi x_2)}{1-\nu^2} \;, \\
\sigma_{22}(x_1,x_2) =& 2 E \pi \frac{ \sin(2 \pi x_1) \cos(2 \pi x_2) + \nu \cos(2 \pi x_1) \sin(2 \pi x_2)}{1-\nu^2} \;, \\
\sigma_{12}(x_1,x_2) =& E \pi \frac{ \cos(2 \pi x_1) \sin(2 \pi x_2) + \sin(2 \pi x_1) \cos(2 \pi x_2)}{1+\nu} \;, \\
f_{1}(x_1,x_2) = f_{2}(x_1,x_2) =& -2 E \pi^2 \frac{ (\nu+1)\cos(2 \pi x_1) \cos(2 \pi x_2) + (\nu-3) \sin(2 \pi x_1) \sin(2 \pi x_2)}{1-\nu^2} \;,
\end{align}
\end{subequations}
on the domain $\Omega \in [-1,1]^2$ with the mapping given by \eqref{eq:mapping}. The particular stress field is chosen to be
\begin{subequations}
\begin{align*}
\sigma_{11}^p =& E \pi \frac{ (\nu+1)\sin(2 \pi x_1) \cos(2 \pi x_2) - (\nu-3) \cos(2 \pi x_1) \sin(2 \pi x_2)}{1-\nu^2} \;,\\
\sigma_{22}^p =& E \pi \frac{ (\nu+1)\cos(2 \pi x_1) \sin(2 \pi x_2) - (\nu-3) \sin(2 \pi x_1) \cos(2 \pi x_2)}{1-\nu^2} \;,\\
\sigma_{12}^p =& 0 \;.
\end{align*}
\end{subequations}

The complementary strain energy, \eqref{eq:comp_strain}, applied to the manufactured solution in \eqref{eq:comp_energy_man_sol} gives the exact value $U_C^{ex} = 58.566883$. In Tables~\ref{tab:strain_energy_stress}, \ref{tab:strain_energy_stress_c0_15} and \ref{tab:strain_energy_stress_c0_30} the complementary strain energy of the approximated solution is shown for $c = \left\lbrace 0,0.15,0.3 \right\rbrace$, and it is observed that the approximated complementary strain energy approaches the exact value from above as expected no matter how distorted the grid is.

\MIG{Note that the convergence of the complementary strain energy is not monotonic.} We expect that when we reduce $h_{el}$ or increase $N$ that the computed strain energy will decrease. As shown in the Tables~\ref{tab:strain_energy_stress}, \ref{tab:strain_energy_stress_c0_15} and \ref{tab:strain_energy_stress_c0_30} this is not always the case. The computed energies which behave contrary to their expected convergence are underlined in these tables.
\begin{table}
\caption{Values of $U_C^h$ in \eqref{eq:comp_strain} using the mapping in \eqref{eq:mapping} with $c=0.0$ for $N=\left\lbrace 2,5,10 \right\rbrace$ and different undeformed element sizes. Forces are specified on all boundaries. Between brackets the value of $U_C^h$ is given when rotation is discretized in the GLL-points as given in \eqref{eq:expan_rot_GLL}.}
\centering
\tabsize
\begin{tabular}{cccc}
\toprule
$h_{el}$ & $N=2$ & $N=5$ & $N=10$ \\
\midrule
$2$ & $81.926894$ $(81.926894)$ & $\underline{82.363976}$ $(77.050078)$ & $58.571086$ $(58.579903)$ \\
$1$ & $81.926894$ $(81.926894)$ & $58.843215$ $(59.374633)$ & $58.566883$ $(58.566884)$ \\
$0.5$ & $64.629318$ $(77.853708)$ & $58.566917$ $(58.570155)$ & $58.566883$ $(58.566883)$ \\
$0.25$ & $58.832253$ $(61.011013)$ & $58.566884$ $(58.566900)$ & $58.566883$ $(58.566883)$\\
$0.125$ & $58.581907$ $(59.166764)$ & $58.566883$ $(58.566883)$ & $58.566883$ $(58.566883)$ \\
\bottomrule
\end{tabular}
\label{tab:strain_energy_stress}
\end{table}
\begin{table}
\caption{Values of $U_C^h$ in \eqref{eq:comp_strain} using the mapping in \eqref{eq:mapping} with $c=0.15$ for $N=\left\lbrace 2,5,10 \right\rbrace$ and different undeformed element sizes. Forces are specified on all boundaries. Between brackets the value of $U_C^h$ is given when rotation is discretized in the GLL-points as given in \eqref{eq:expan_rot_GLL}.}
\centering
\tabsize
\begin{tabular}{cccc}
\toprule
$h_{el}$ & $N=2$ & $N=5$ & $N=10$ \\
\midrule
$2$ & $81.926894$ $(81.926894)$ & $\underline{100.127567}$ $(63.185431)$ & $65.932457$ $(62.042183)$ \\
$1$ & $\underline{98.520452}$ $(64.683021)$ & $62.336805$ $(\underline{67.838096})$ & $58.569756$  $(58.601725)$ \\
$0.5$ & $73.321298$ $(\underline{76.952091})$ & $58.617041$ $(59.083907)$ & $58.566883$ $(58.566887)$ \\
$0.25$ & $59.961965$ $(\underline{78.325747})$ & $58.566974$ $(58.569018)$ & $58.566883$ $(58.566883)$ \\
$0.125$ & $58.651790$ $(\underline{74.791834})$ & $58.566883$ $(58.566902)$ & $58.566883$ $(48.566883)$ \\
\bottomrule
\end{tabular}
\label{tab:strain_energy_stress_c0_15}
\end{table}
\begin{table}
\caption{Values of $U_C^h$ in \eqref{eq:comp_strain} using the mapping in \eqref{eq:mapping} with $c=0.3$ for $N=\left\lbrace 2,5,10 \right\rbrace$ and different undeformed element sizes. Forces are specified on all boundaries. Between brackets the value of $U_C^h$ is given when rotation is discretized in the GLL-points as given in \eqref{eq:expan_rot_GLL}.}
\centering
\tabsize
\begin{tabular}{cccc}
\toprule
$h_{el}$ & $N=2$ & $N=5$ & $N=10$ \\
\midrule
$2$ & $81.926894$ $(81.926894)$ & $\underline{114.733969}$ $(68.884233)$ & $\underline{83.029888}$ $(63.403240)$ \\
$1$ & $\underline{89.196288}$ $(76.033509)$ & $78.835567$ $(\underline{71.435962})$ & $58.800555$ $(59.033822)$ \\
$0.5$ & $\underline{87.627401}$ $(\underline{78.637465})$ & $59.397067$ $(61.098091)$ & $58.566907$ $(58.567109)$ \\
$0.25$ & $65.070084$ $(\underline{80.383044})$ & $58.570680$ $(58.611592)$ & $58,566883$ $(58.566883)$ \\
$0.125$ & $59.066273$ $(\underline{76.980788})$ & $58.566887$ $(58.567241)$ & $58.566883$ $(58.566883)$ \\
\bottomrule
\end{tabular}
\label{tab:strain_energy_stress_c0_30}
\end{table}
\begin{figure}
\centering
\includegraphics[width=1.00\textwidth]{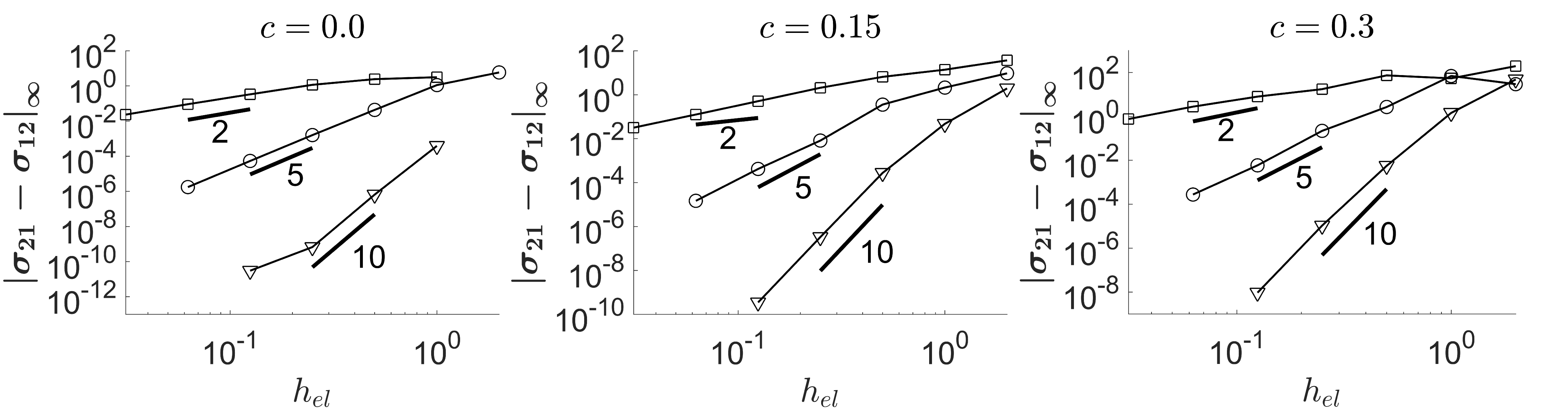} %
\caption{Convergence of $\|\sigma_{12}-\sigma_{21}\|_{\infty}$ with respect to the undeformed element size $h_{el}$ using the mapping in \eqref{eq:mapping} with $c=0.15$. $\Box:N=2$, $\circ:N=5$, $\bigtriangledown:N=10$.}
\label{fig:convergence_symmetry_all_grids}
\end{figure}
This irregular convergence behavior may be attributed to the fact that rotational equilibrium, \eqref{eq:cons_ang_mom} is only weakly imposed, see Figure~\ref{fig:convergence_symmetry_all_grids}. By discretizing the rotation field in the Gauss points we do not have enough Lagrange multipliers $\omega$ in \eqref{eq:Elastic_lagrangian} to strongly enforce \KRO{equilibrium of moments}.
The rotation field converges to the exact solution with mesh refinement for both $c=0.0$, Figure~\ref{fig:convergence_rotation_c0} and on the curvilinear mesh with $c=0.15$, Figure~\ref{fig:convergence_rotation_c15}.
\begin{figure}
\centering
\includegraphics[width=1.00\textwidth]{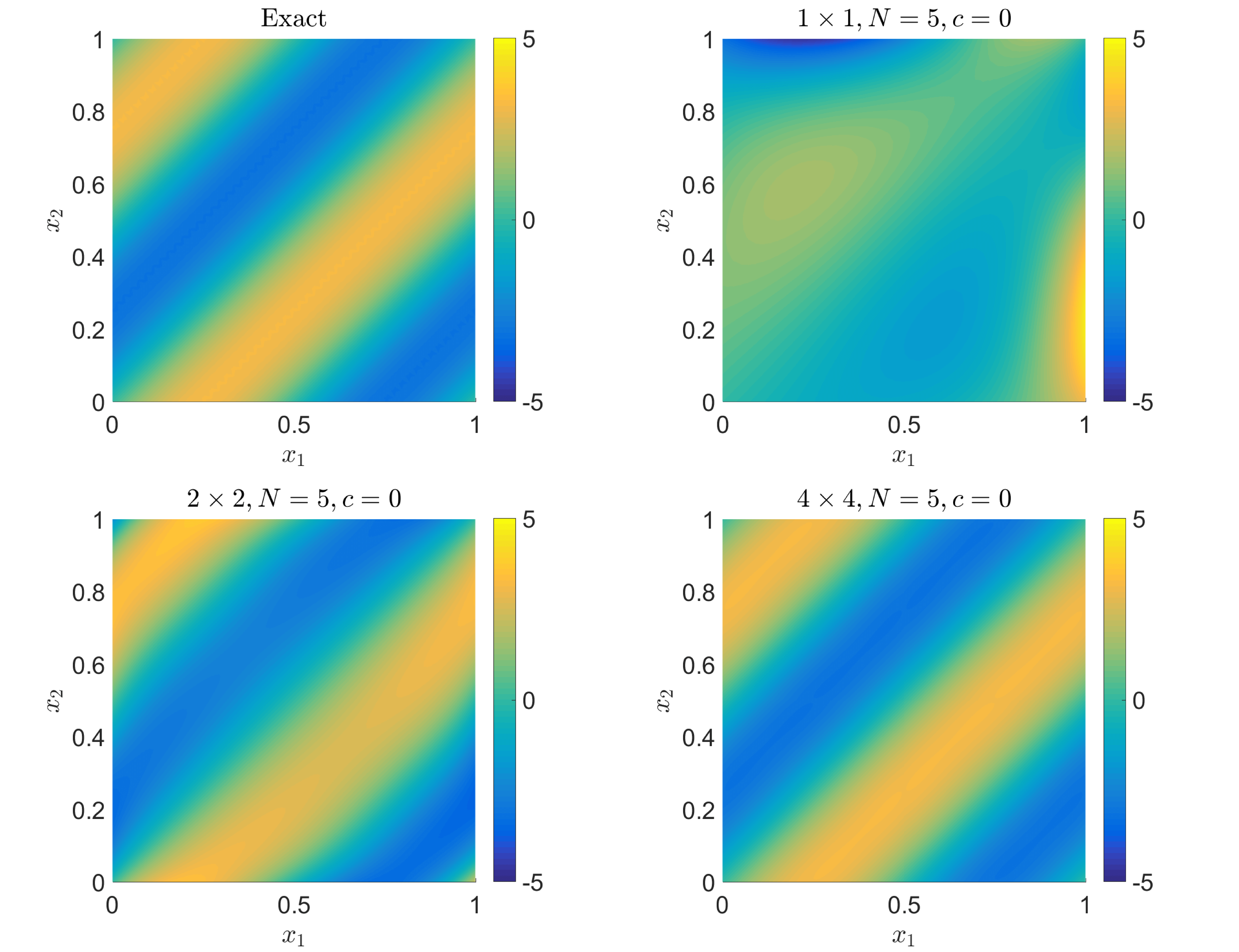} %
\caption{Convergence of the rotation field for $N=5$ on an orthogonal grid, $c=0.0$. The exact solution (top left), approximate solution in 1 element (top right), approximate solution in $2 \times 2$ elements (bottom left) and approximate solution in $4 \times 4$ elements (bottom right).}
\label{fig:convergence_rotation_c0}
\end{figure}
\begin{figure}
\centering
\includegraphics[width=1.00\textwidth]{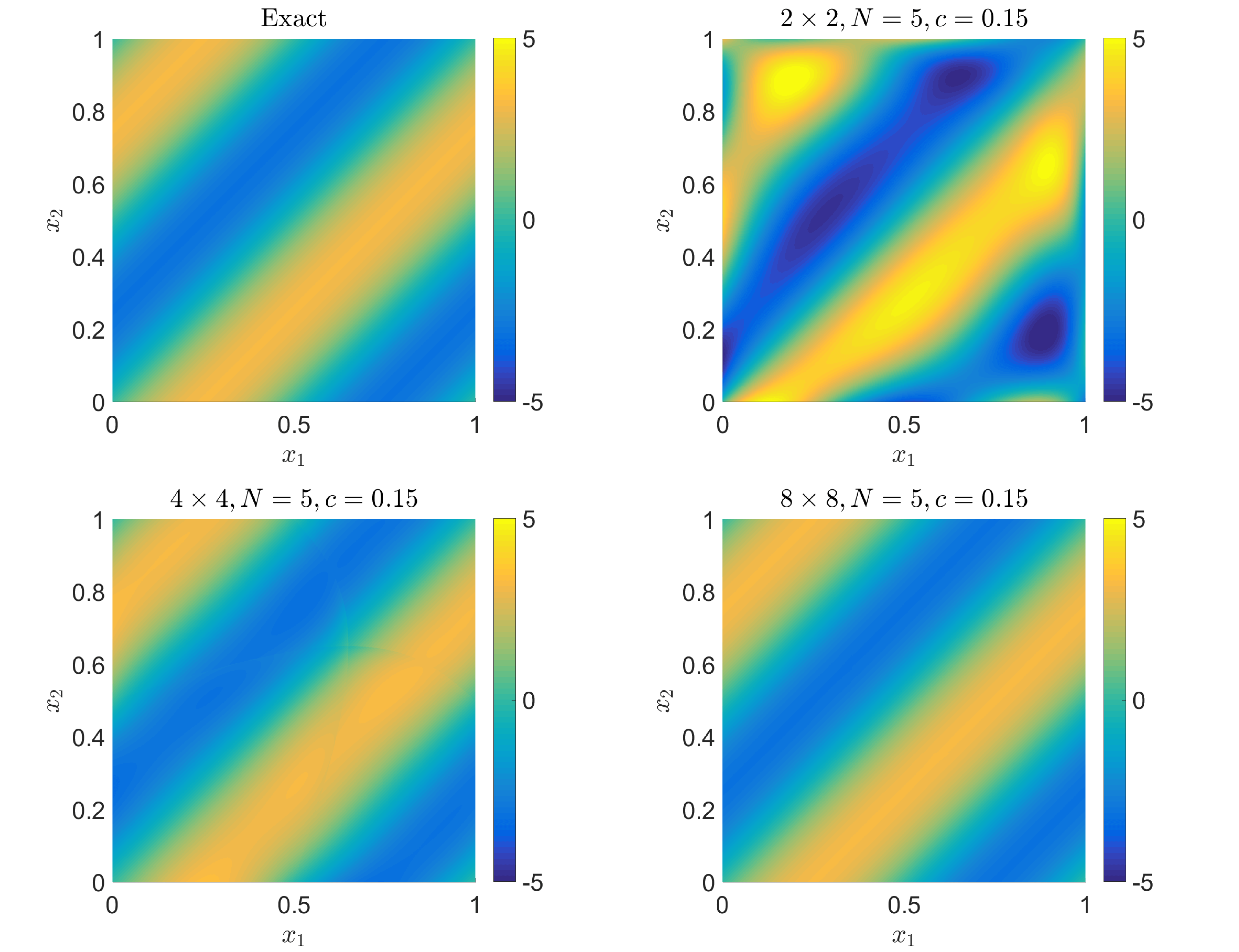} %
\caption{Convergence of the rotation field for $N=5$ on a curvilinear grid, $c=0.15$. The exact solution (top left), approximate solution in $2 \times 2$ elements (top right), approximate solution in $4 \times 4$ elements (bottom left) and approximate solution in $8 \times 8$ elements (bottom right).}
\label{fig:convergence_rotation_c15}
\end{figure}

If we discretize the rotations in the GLL points, that is,
\begin{subequations} \label{eq:expan_rot_GLL}
\begin{align}
\omega_{1}^{s,h}(\xi_1,\xi_2,\xi_3) =& \sum\limits_{i=0}^{N} \sum\limits_{j=0}^{N} \sum\limits_{k=0}^{N} \left( O_{1}^s \right)_{i,j,k} h_i(\xi_1) h_j(\xi_2) h_k(\xi_3) \;,
\\
\omega_{2}^{s,h}(\xi_1,\xi_2,\xi_3) =& \sum\limits_{i=0}^{N} \sum\limits_{j=0}^{N} \sum\limits_{k=0}^{N} \left( O_{2}^s \right)_{i,j,k} h_i(\xi_1) h_j(\xi_2) h_k(\xi_3) \;,
\\
\omega_{3}^{s,h}(\xi_1,\xi_2,\xi_3) =& \sum\limits_{i=0}^{N} \sum\limits_{j=0}^{N} \sum\limits_{k=0}^{N} \left(O_{3}^s \right)_{i,j,k} h_i(\xi_1) h_j(\xi_2) h_k(\xi_3) \;,
\end{align}
\end{subequations}
where
\[ \left( O_{1}^s \right)_{i,j,k} = \omega_{1}^{s,h}(\left( \xi_1 \right)_{2},\left( \xi_2 \right)_{j},\left( \xi_3 \right)_{k}) \;,\;\;\; \left( O_{2}^s \right)_{i,j,k} = \omega_{2}^{s,h}(\left( \xi_1 \right)_{2},\left( \xi_2 \right)_{j},\left( \xi_3 \right)_{k}) \;,\]
\[ \left( O_{3}^s \right)_{i,j,k} = \omega_{3}^{s,h}(\left( \xi_1 \right)_{2},\left( \xi_2 \right)_{j},\left( \xi_3 \right)_{k}) \;,\]
then rotational equilibrium is more strongly weighted in the variational formulation. The complementary strain energy obtained by using the rotation in the GLL points are listed in the Tables~\ref{tab:strain_energy_stress}, \ref{tab:strain_energy_stress_c0_15} and \ref{tab:strain_energy_stress_c0_30} between brackets. Once again the computed energies which do not converge monotonically are underlined in these tables. We see that in Table~\ref{tab:strain_energy_stress}, the results for the undistorted mesh, all computed complementary strain energies converge monotonically from above to the exact strain energy. Table~\ref{tab:strain_energy_stress_c0_15} and Table~\ref{tab:strain_energy_stress_c0_30}  reveal that on the deformed grids monotonic converge no longer takes place.
\begin{figure}
\centering
\includegraphics[width=1.00\textwidth]{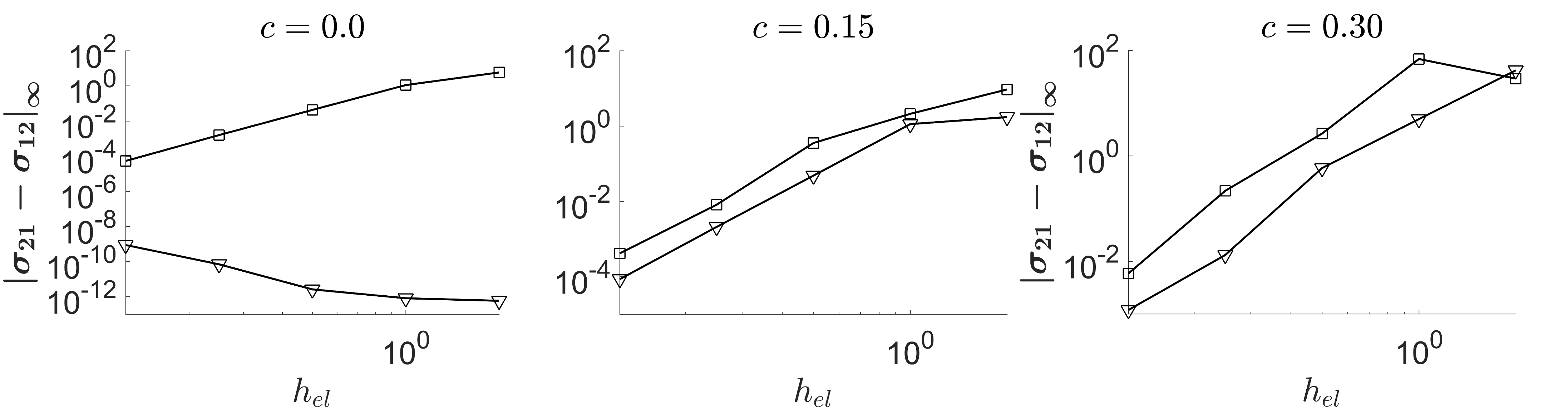} %
\caption{Convergence of $\|\sigma_{12}-\sigma_{21}\|_{\infty}$ for $N=5$ with respect to the undeformed element size $h_{el}$ using the mapping in \eqref{eq:mapping} with $c=0.15$. $\Box:$ rotation in the Gauss points, \eqref{eq:expan_rot} $\bigtriangledown:$ rotation in the Gauss-Lobatto points, \eqref{eq:expan_rot_GLL}.}
\label{fig:Comparison_G_GL_points_for_rot}
\end{figure}
Figure~\ref{fig:Comparison_G_GL_points_for_rot} shows that for $c=0.0$ symmetry of the stress tensor is strongly imposed, which agrees with monotonic convergence from above towards the analytic strain energy. Figure~\ref{fig:Comparison_G_GL_points_for_rot} shows that on the deformed grid, $c=0.15$ and $c=0.30$, symmetry of the stress tensor is only weakly enforced. On the deformed grids monotonic convergence towards the analytic strain energy is lost; see also Figure~\ref{fig:d_s_str}.

However, it is observed that the matrix system becomes singular and the rank deficiency increases as the number of elements increases, when the rotation is discretized in the GLL-points. The appearance of singular kinematic modes in equilibrium models is well known, see for instance \cite[Ch.6]{bookPian} and \cite[Ch.5\&6]{BookAlmeidaMaunder}.
Further investigations of a proper expansion of the rotation field is required in future work to get equilibrium of moments on deformed elements and no rank deficiency in the matrix system.
\begin{figure}
\centering
\includegraphics[width=1.00\textwidth]{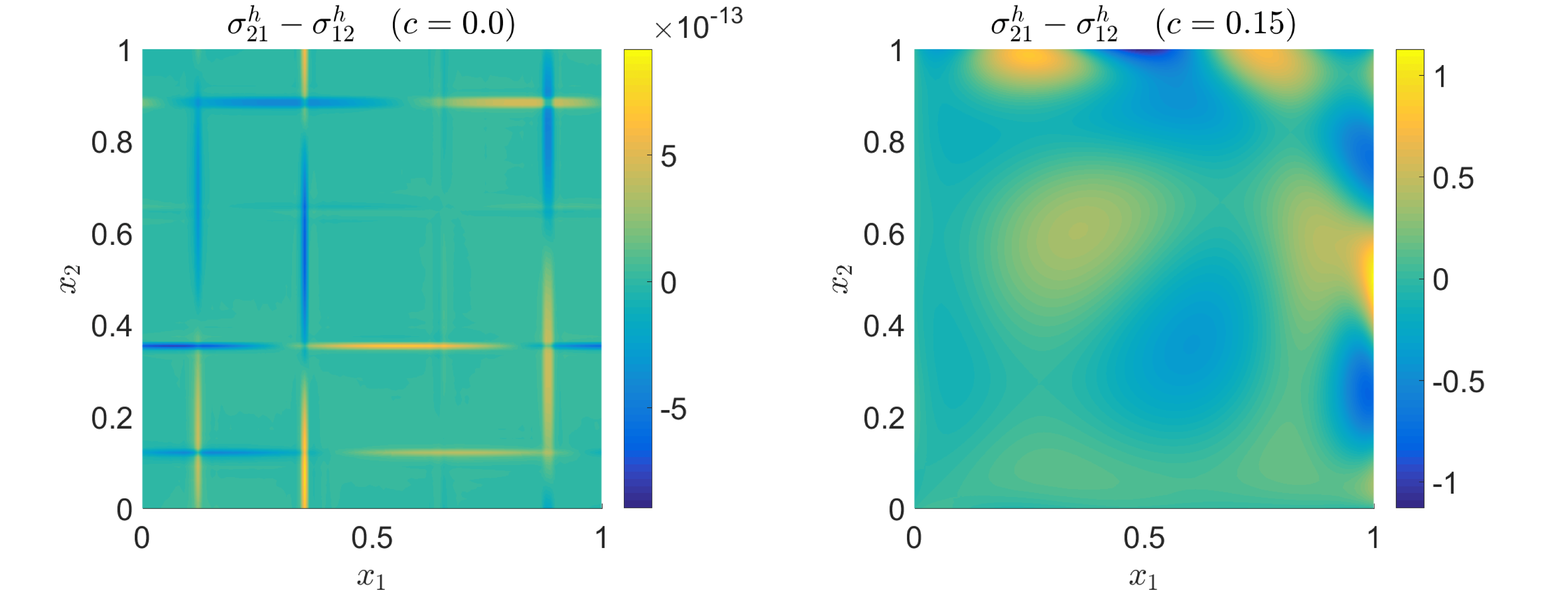}
\caption{Symmetry of the stress tensor for $c=0.0$ to the left and for $c=0.15$ on the right on a $2\times 2$ grid with polynomial approximation $N=5$.}
\label{fig:d_s_str}
\end{figure}
In Figure~\ref{fig:d_s_str} the polynomial solutions are sampled in a large number of points to produce the plots.
\begin{figure}
\centering
\includegraphics[width=0.50\textwidth]{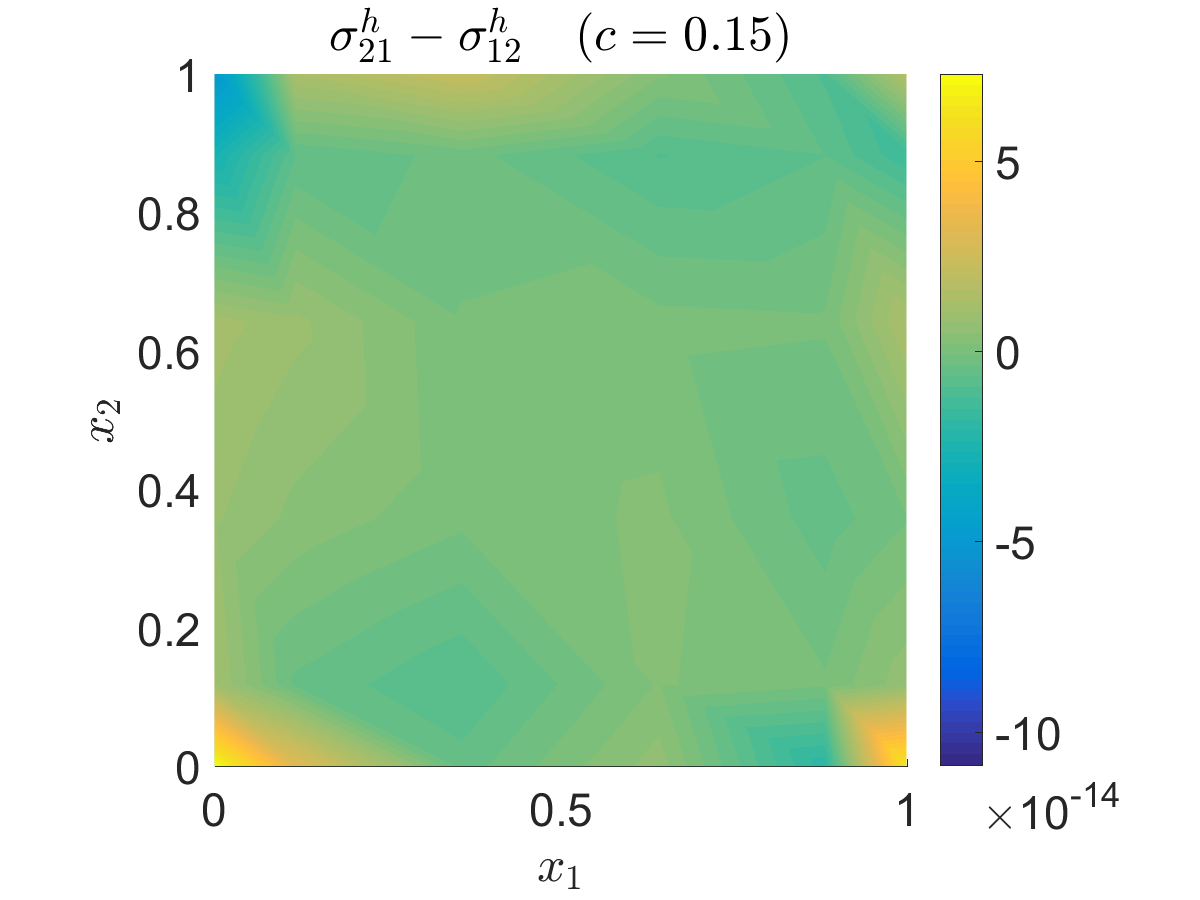}
\caption{Symmetry of the stress tensor for $c=0.15$ on a $2\times 2$ grid with polynomial approximation $N=5$ sampled in the Gauss-Lobatto-Legendre points.}
\label{fig:d_s_str2}
\end{figure}
If we only sample the solution in the Gauss-Lobatto points, then the contour plot for $\sigma^h_{21}-\sigma^h_{12}$ is displayed in Figure~\ref{fig:d_s_str2}. This figure illustrates that symmetry of the stress tensor is satisfied up to $O(10^{-14})$ in the integration points. These two plots demonstrate that symmetry of the stress tensor is restored pointwsie on orthogonal grids and only in the integration points on a curvilinear grid. The price one pays for this restored symmetry is that the system matrix is singular and the rotation field is non-unique.

\section{Comparison with Traditional Displacement FE Method} \label{sec:FE_comp}
\begin{figure}
\centering
\includegraphics[width=1.00\textwidth]{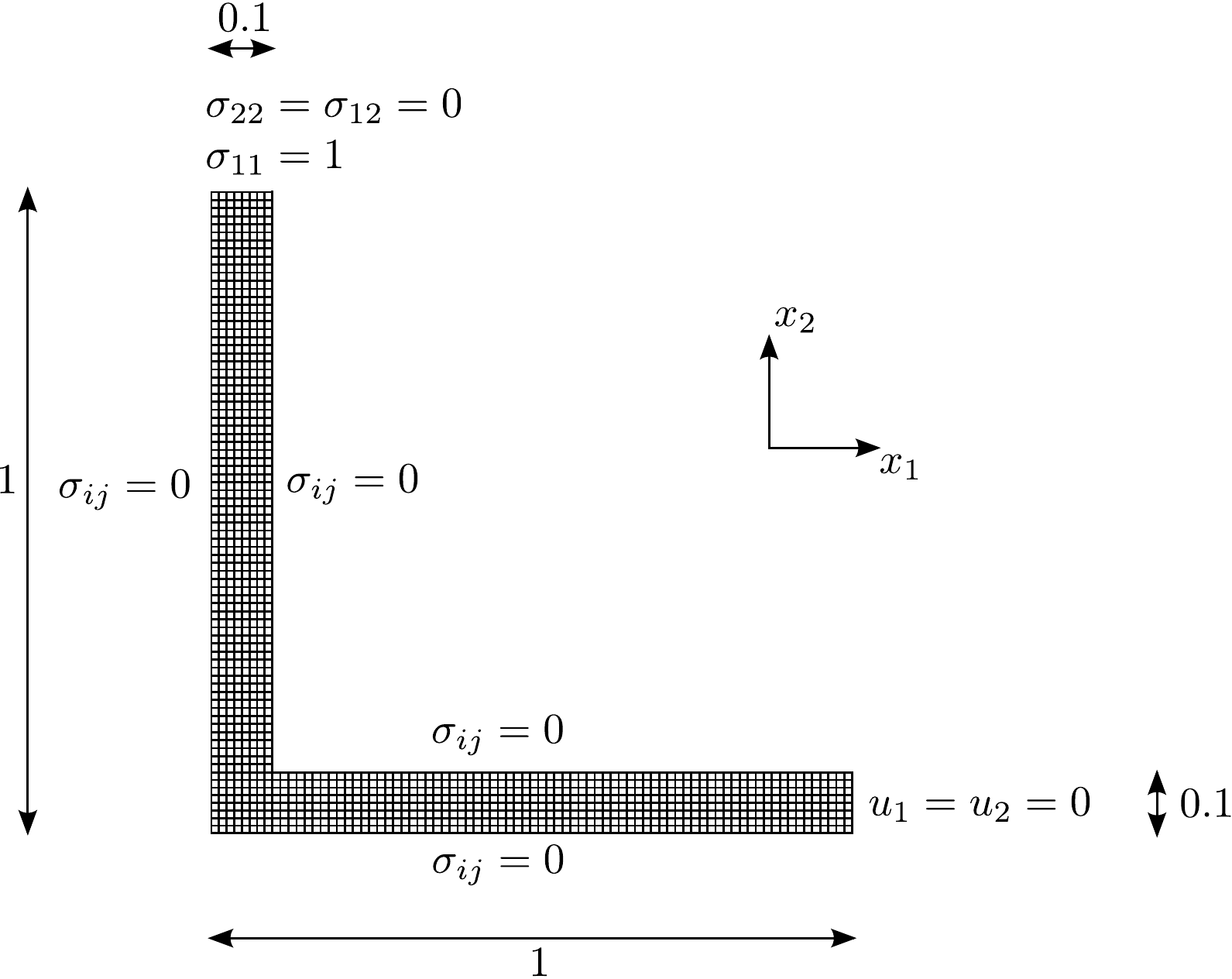}
\caption{L shape domain with boundary conditions.}
\label{fig:L_shape}
\end{figure}

\KRO{In this section the method is compared to a traditional displacement based Galerkin FE method, see for instance \cite[Ch. 3.3]{Cook_2001}. The comparison is done by considering the solution on an L-shaped domain with the size, grid and boundary conditions shown in Figure~\ref{fig:L_shape}. This test case has a singularity of the stress components at the intersection of the legs, and it will therefore show how the presented method behaves near singularities. The material properties are set to $E=1$ and $\nu=0.3$, and the elements are quadratic with a size of 0.0125, which produces 1216 elements.}

\begin{figure}
\centering
\includegraphics[width=1.0\textwidth]{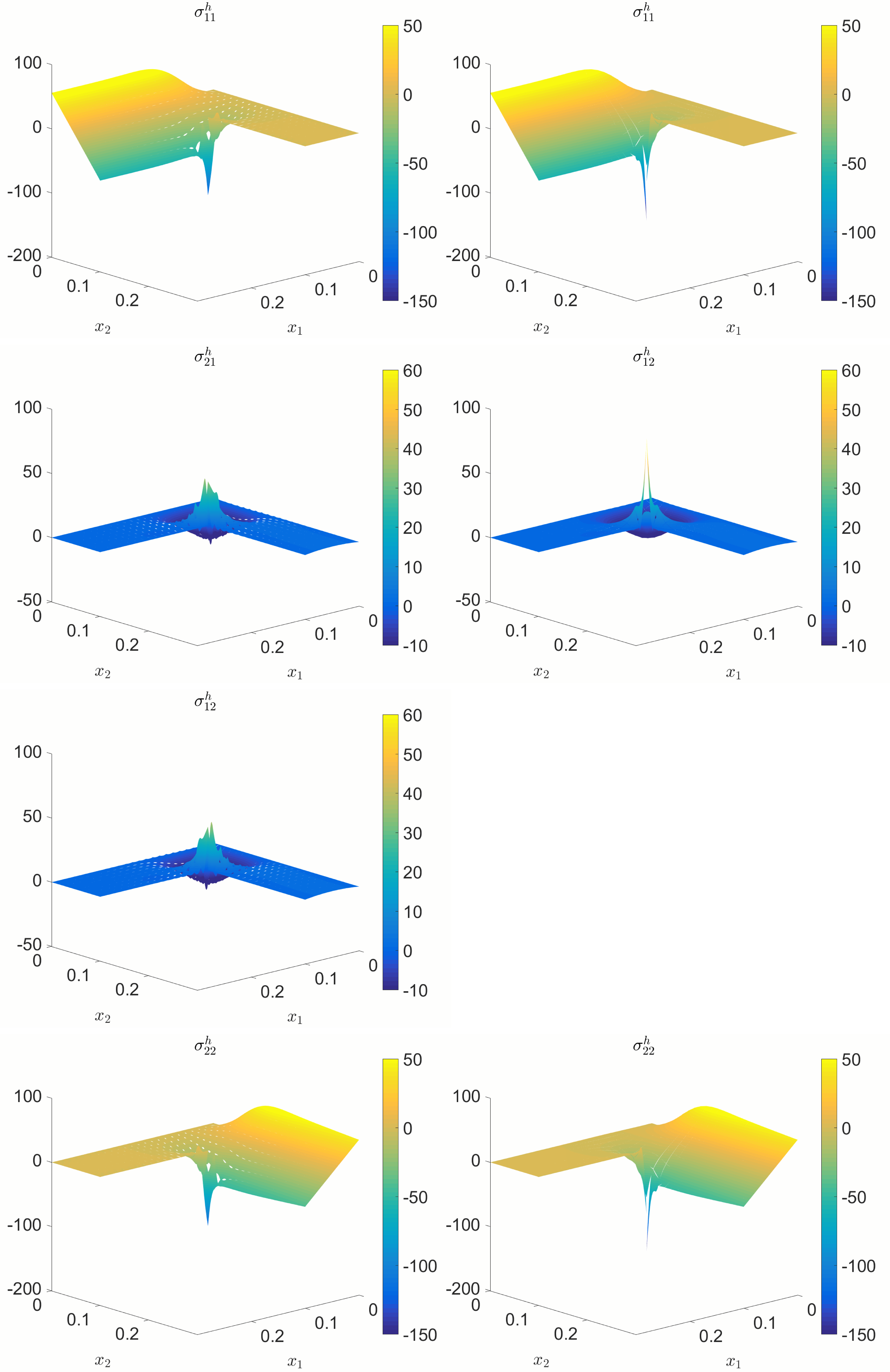}
\caption{The stress components at the singularity of the L-shape test case for $N=2$. Left column: The method presented in the paper. Right column: Traditional displacement based FE method with Q4 elements.}
\label{fig:L_shape_N2}
\end{figure}

\begin{figure}
\centering
\includegraphics[width=1.0\textwidth]{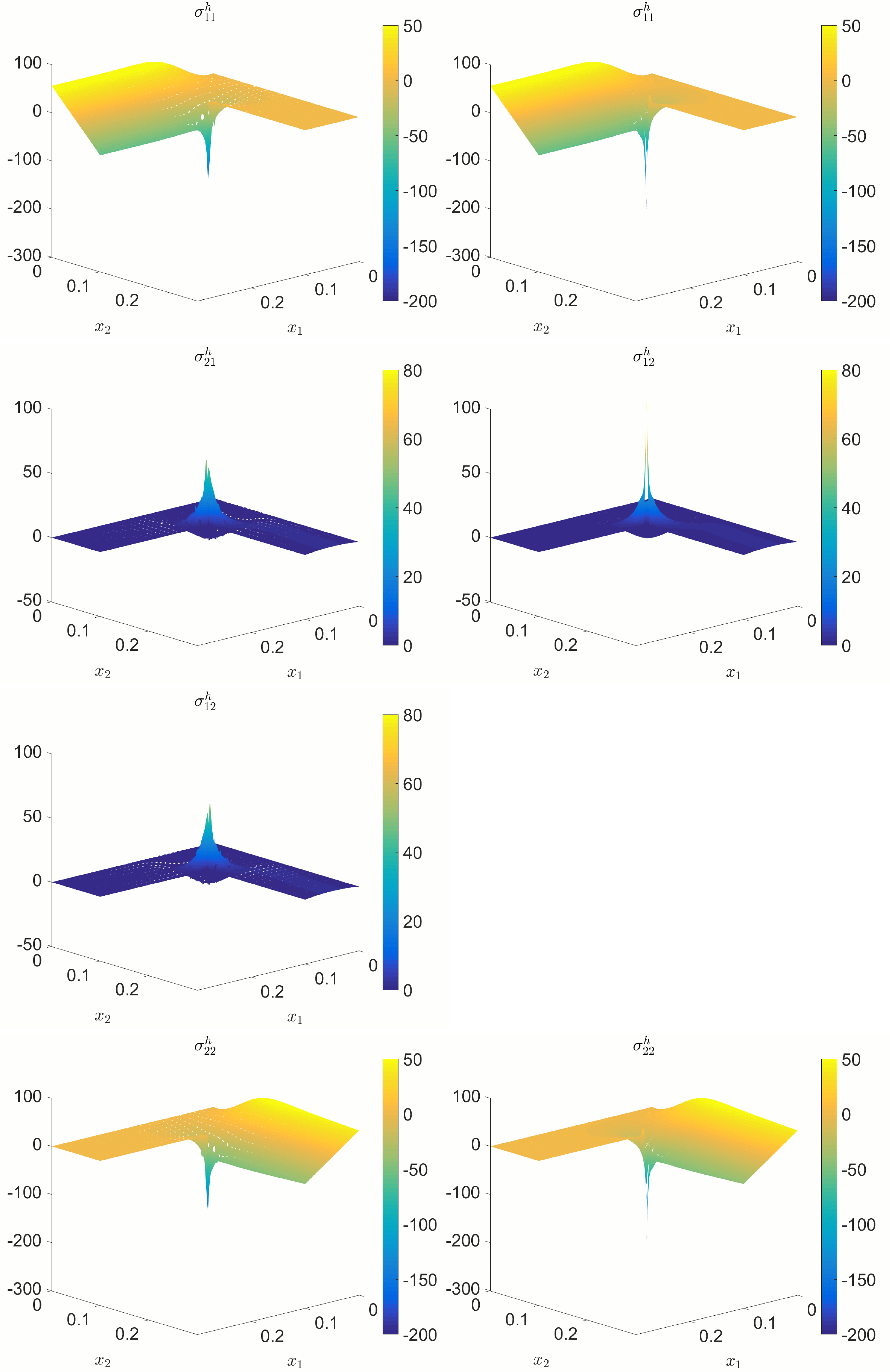}
\caption{The stress components at the singularity of the L-shape test case for $N=3$. Left column: The method presented in the paper. Right column: Traditional displacement based FE method with Q9 elements.}
\label{fig:L_shape_N3}
\end{figure}

\KRO{The stress components for $N=2$ are plotted in Figure~\ref{fig:L_shape_N2}, where the results for the current formulation are in the left column and the FE method in the right column. In the FE method $N=2$ corresponds to Q4 elements. A similar plot is given in Figure~\ref{fig:L_shape_N3}, which in a traditional FE method corresponds to Q9 elements. In both plots the singularity is clearly observed. The magnitude of the stress components in the current formulation is smaller than in the FE method. The reason for this is that in the FE method a computational node is located at the singularity, while in the current formulation the DOFs for the stress components are the force components located on the boundaries of the element. This property also implies that the stress components in the current formulation is continuous in one direction, but not in the other direction. It is furthermore observed that in the FE method the shear stress component displays two spikes at the singularity, whereas the current formulation only shows one peak.}

\begin{figure}
\centering
\includegraphics[width=1.0\textwidth]{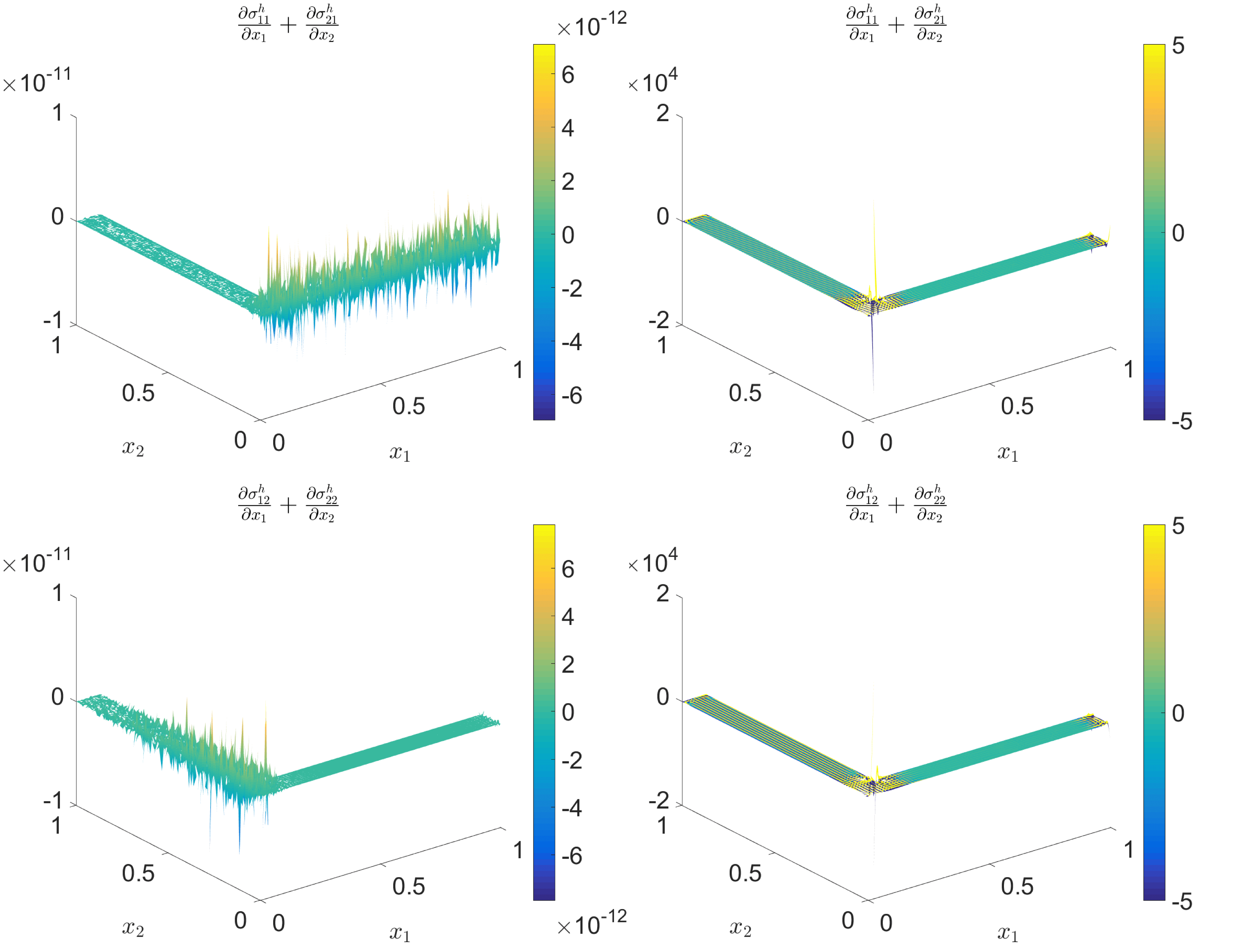}
\caption{The force equilibrium of the L-shape test case for $N=2$. Left column: The method presented in the paper. Right column: Traditional displacement based FE method with Q4 elements.}
\label{fig:L_F_eq_shape_N2}
\end{figure}

\begin{figure}
\centering
\includegraphics[width=1.0\textwidth]{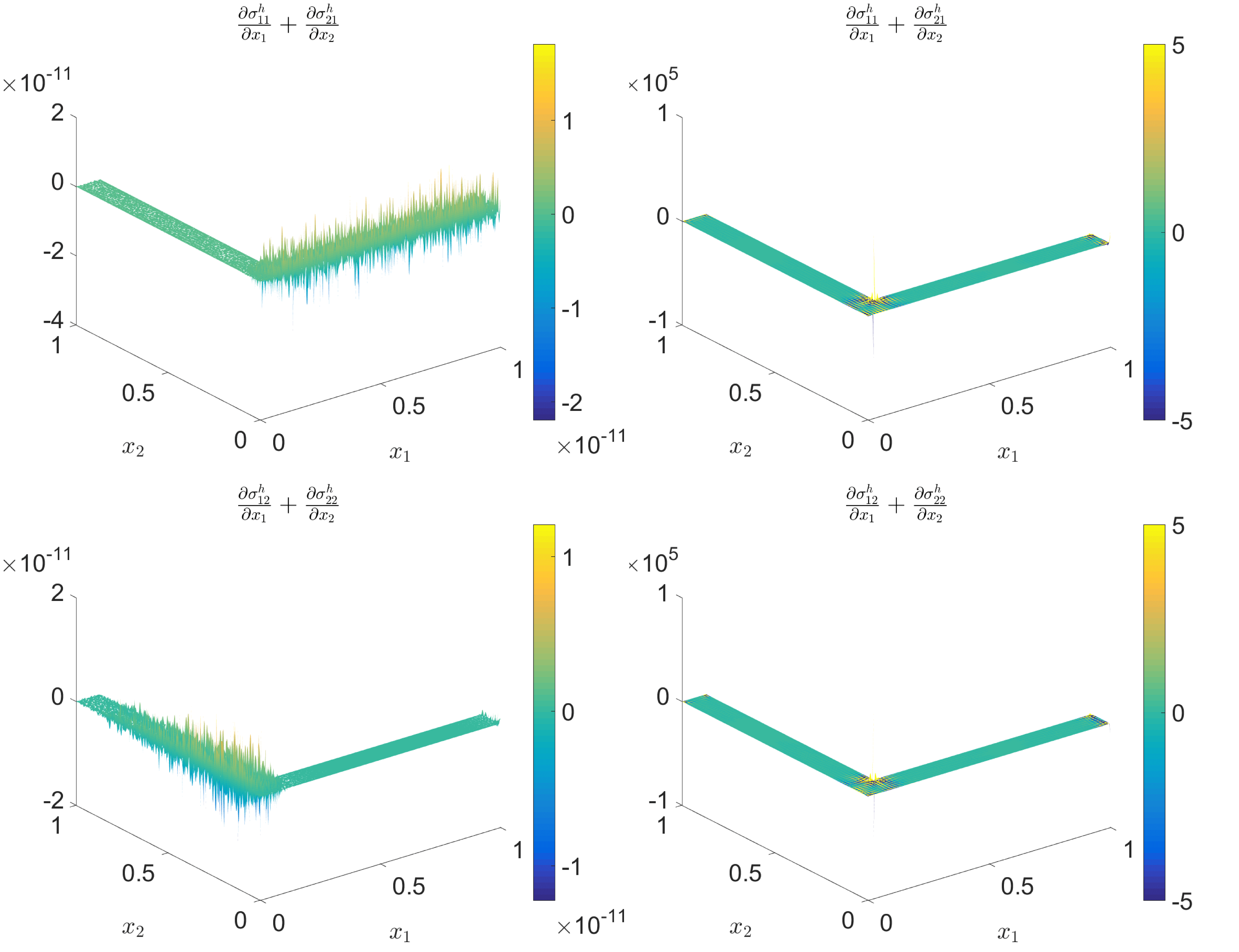}
\caption{The force equilibrium of the L-shape test case for $N=3$. Left column: The method presented in the paper. Right column: Traditional displacement based FE method with Q9 elements.}
\label{fig:L_F_eq_shape_N3}
\end{figure}

\KRO{In Figure~\ref{fig:L_F_eq_shape_N2} and Figure~\ref{fig:L_F_eq_shape_N3} the force equilibrium equations are evaluated over the domain for $N=2$ and $N=3$, respectively. It is clearly seen that the current formulation satisfies the force equilibrium equations to machine precision both for $N=2$ and $N=3$, while the FE method has residuals between -5 and 5 in the majority of the domain and at the singularity they are between $10^4$ and $10^5$, which are observed by the thin spikes at the singularity.}

\begin{figure}
\centering
\includegraphics[width=1.00\textwidth]{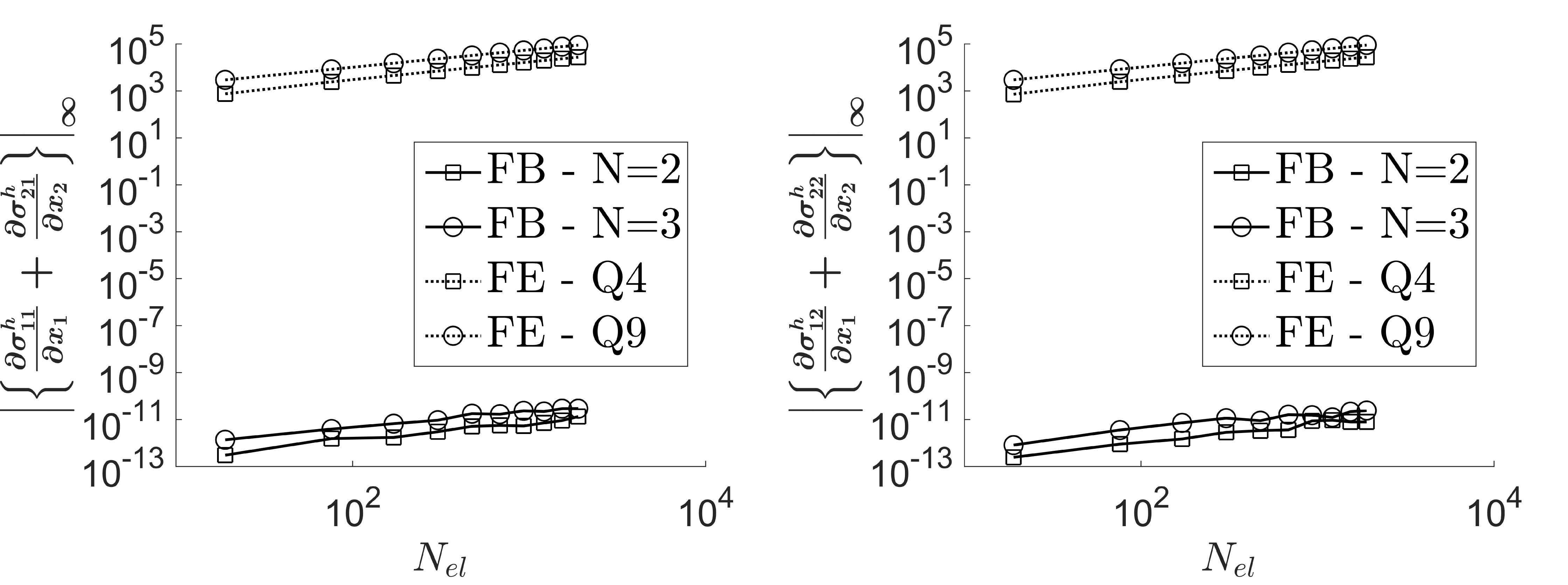}
\caption{Comparison of the force equilibrium between the method presented in this paper (FB) and a traditional displacement FE method (FE).}
\label{fig:L_shape_force_eq}
\end{figure}

\KRO{Next, the equations are solved on different mesh sizes ranging from element sizes of 0.1 to 0.0125. In Figure~\ref{fig:L_shape_force_eq} the maximum inequality of the equilibrium equations is plotted as a function of the number of elements. As observed the singularity does not have an influence on the equilibrium of translational forces even for large element sizes. The plots indicate the maximum value of the equilibrium equations evaluated in $100 \times 100$ points in each element. The traditional FE method have an inequality in the order of $10^3$ to $10^5$.}

\begin{figure}
\centering
\includegraphics[width=0.50\textwidth]{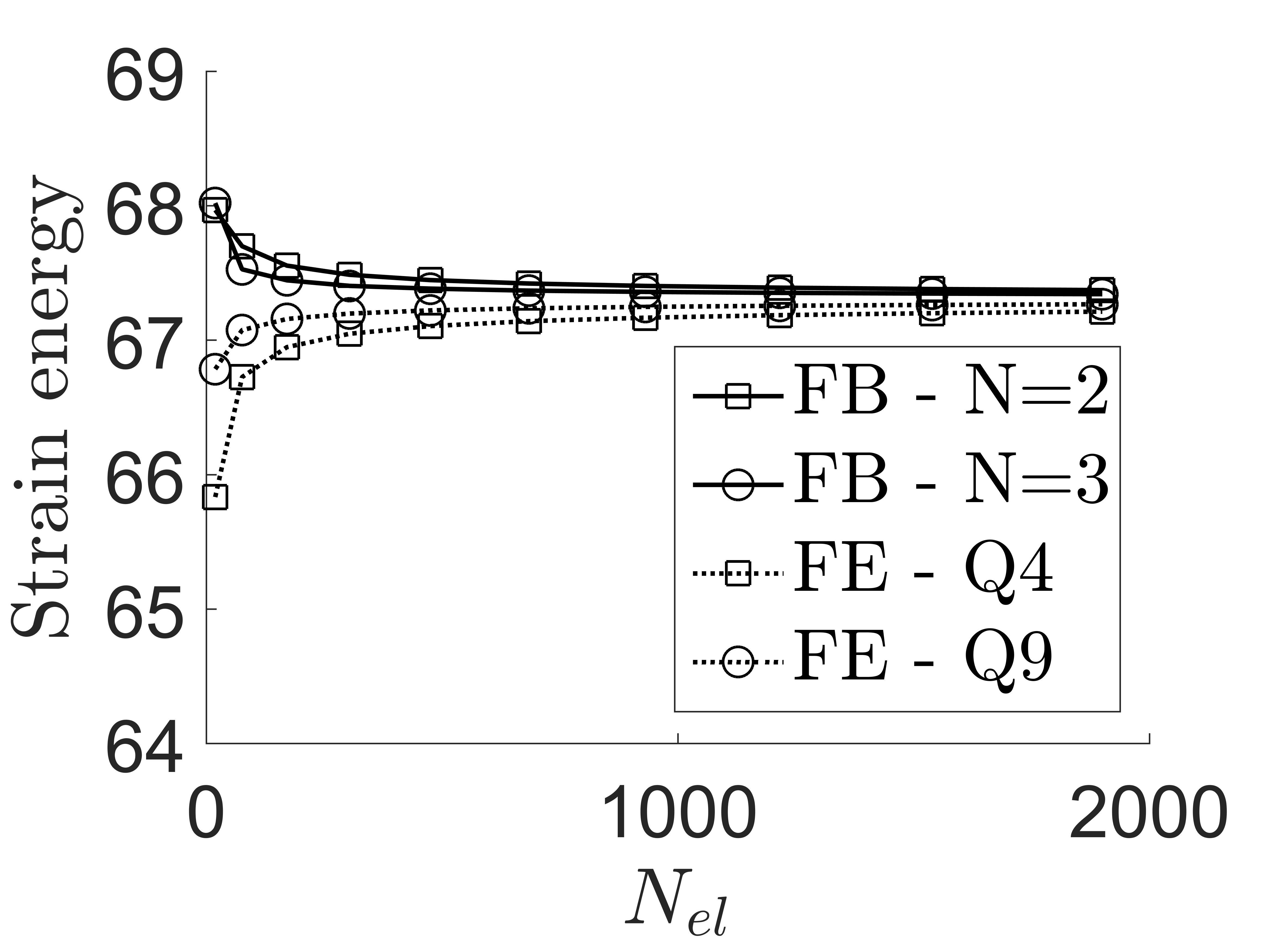}
\caption{Convergence of the strain energy of the method presented in this paper (FB) and a traditional displacement FE method (FE).}
\label{fig:convergence_L_shape_se}
\end{figure}

\KRO{As discussed in Section~\ref{sec:Comp_energy} the force based methods will cause the complementary strain energy to converge from above. Displacement based methods will, on the other hand, have the strain energy converge from below, see for instance \cite{Veubeke_1964,Kempeneers_2010}. These trends are clearly observed in Figure~\ref{fig:convergence_L_shape_se}, where the strain energy for both methods is plotted as a function of the number of elements.}

\begin{figure}
\centering
\includegraphics[width=1.00\textwidth]{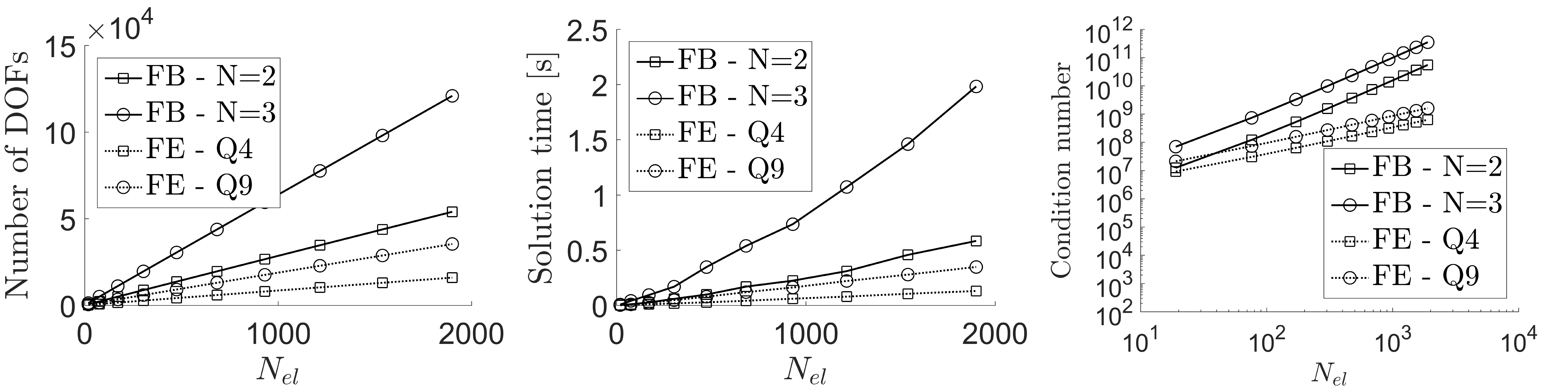}
\caption{Comparison of number of DOFs (left), solution time (middle) and condition number (right) between the method presented in this paper (FB) and a traditional displacement FE method (FE).}
\label{fig:L_shape_sol_time}
\end{figure}

\KRO{The current formulation presented in this paper intoduces additional DOFs compared to the traditional FE method, see Figure~\ref{fig:L_shape_sol_time}. In the tradtional FE method the displacement components are solved for, while the method in this paper also solves for the surface force components of the elements and rotation values. This is clearly seen in Figure~\ref{fig:L_shape_sol_time}, where the solution time and the condition number of the system matrix are plotted as a function of the number of elements. It is clearly observed that the current formulation has both longer solution time and larger condition numbers than the FE method. It is also observed that the solution time for the current formulation increases considerably as $N$ increases. In \cite{Jain} the sparsity of the equation system is considerably increased, and thereby reducing the condition number, through the use of dual polynomial spaces. Implementing such polynomials in the method presented in this method should reduce the solution time and condition numbers. A saddle point decomposition, like the Uzawa algorithm, see \cite{Karniadakis_2005}, could also be investigated. Here the system of equations is broken down into a number of smaller Poisson equations, which is easier to solve. These methods are not investigated in the present paper, but are addressed in future work.}

\section{Conclusion} \label{sec:conclusions}
In this paper a spectral element method is presented, which satisfies the translational equilibrium of forces in a structural problem to machine precision, independent of the size and shape of the individual elements, for piecewise polynomial body forces. Optimal convergence rates of the problem based on the order of the polynomial expansions are obtained on orthogonal grid and curvilinear grids. Through the use of \emph{edge expansion polynomials} discrete surface force components and body force components acting on the elements are connected to the stress components and body force density components. The discrete representation of the force equilibrium equations is just a summation of the force components, and yields exact force equilibrium both internally and between the elements when no body forces are present.

When the rotation is discretized in the Gauss points, the stress tensor is only weakly symmetric. Convergence of the complementary strain energy is from above towards the exact energy, but this convergence is non-monotonic. In this case the rotation field converges to the exact solution on orthogonal, Figure~\ref{fig:convergence_rotation_c0}.

When the rotation is discretized in the Gauss-Lobatto points, the stress tensor is strongly symmetric on an orthogonal domain and the monotonic convergence from above towards the analytic complementary strain energy is restored for the analytic test problem considered in this paper. Uniqueness of the discrete solution is lost due to the presence of singular modes in the rotation field.

On curvilinear grids, see Figure~\ref{fig:mesh_3}, both discretization of rotation on the Gauss grid and on the Gauss-Lobatto grid lead to weakly symmetric stress tensors, see Figure~\ref{fig:Comparison_G_GL_points_for_rot}. Monotonic convergence of the complementary strain energy is lost, see Tables~\ref{tab:strain_energy_stress_c0_15} and \ref{tab:strain_energy_stress_c0_30}. When a Gauss grid is used for the rotation field, the rotation field converges to the analytic solution, Figure~\ref{fig:convergence_rotation_c15}. On a Gauss-Lobatto grid, symmetry of the stress tensor is only enforced in the integration points, Figure~\ref{fig:d_s_str2}, but outside the integration points the stress tensor is not symmetric, Figure~\ref{fig:d_s_str}.

An extension of the ideas presented in this paper, to a method which will also satisfy pointwise \KRO{equilibrium of moments} without the introduction of spurious kinematic modes will be presented in a follow-up paper.

\vskip 0.3cm
\noindent
\textbf{Acknowledgements} The authors want to thank the reviewers who, by their critical remarks and suggestions, have contributed significantly to the current manuscript.


\begin{thebibliography}{9}

\bibitem{Zienkiewicz_1972} Zienkiewicz~OC. \emph{The Finite Element Method in Engineering Science}.
McGraw-Hill Inc. 1972.

\bibitem{Oden_1976} Oden~JT, Reddy~JN. \emph{Mathematical Theory of Finite Elements}.
John Wiley and sons. 1976.

\bibitem{Bathe_1976} Bathe~KJ, Wilson~EL. \emph{Numerical Methods in Finite Element Analysis}.
Prentice-Hall. 1976.

\bibitem{Reddy_2006} Reddy~JN. \emph{An Introduction to the Finite Element Method} (3rd~edn).
McGraw-Hill. 2006.

\bibitem{Reddy_2010} Reddy~JN, Gartling~DK. \emph{The Finite Element Method in Heat Transfer and Fluid Dynamics}.
CRC Press. 2010.

\bibitem{Monk_2003} Monk~P. \emph{Finite Element Methods for Maxwell's Equations}.
Oxford Science Publications. 2003.

\bibitem{Reddy_2015} Reddy~JN, Srinivasa~A. On the Force-Displacement Characteristics of Finite Elements for Elasticity and Related Problems. \emph{Finite Elements in Analysis and Design 2015}; \textbf{104}:35--40.

\bibitem{Almeida_Continuity} Moitinho~de~Almeida~JP, Texeira~de~Freitas~JA. Continuity conditions for finite element analysis of solids. \emph{International Journal for Numerical Methods in Engineering}; \textbf{33}:845--853, 1992.

\bibitem{Cook_2001} Cook~RD, Malkus~DS, Plesha~ME, Witt~RJ. \emph{Concepts and Applications of Finite Element Analysis} (4th~edn).
John Wiley and sons. 2001.

\bibitem{Marsden_1983} Marsden~JE, Hughes~TJR. \emph{Mathematical Foundations of Elasticity}.
Prentice-Hall. 1983.

\bibitem{Almeida_1996} Moitinho~de~Almeida~JP, Almeida~Pereira~OJB. A Set of Hybrid Equilibrium Finite Element Models for the Analysis of Three-dimensional Solids. \emph{International Journal for Numerical Methods in Engineering 1996}; \textbf{39}:2789--2802.

\bibitem{Veubeke_1964} Fraeijs~de~Veubeke~BM. Upper and Lower Bounds in Matrix Structural Analysis. \emph{AGARDograf 1964}; \textbf{72}:165--201.

\bibitem{Veubeke_1965} Fraeijs~de~Veubeke~BM. Displacements and Equilibrium Models in the Finite Elements Method. \emph{In Stress Analysis, Zienkiewicz OC, Holister GS (eds), Chapter 9}.
Wiley. 1965.

\bibitem{Veubeke_Millard_1976} Fraeijs~de~Veubeke~BM, Millard~A.. Discretization of Stress Fields in the Finite Element Method. \emph{Journal of the Franklin Institue}. Vol. 302 No. 5\&6, pp.~389--412, 1976.

\bibitem{Veubeke_1980} Fraeijs~de~Veubeke~BM. Diffusive Equilibrium Models. \emph{In B.M. Fraeijs de Veubeke memorial volume of selected papers}. Sijthoff \& Noordhoff. 1980.

\bibitem{Almeida_1991} Moitinho~de~Almeida~JP, Freitas~JAT. Alternative Approach to the Formulation of Hybrid Equilibrium Finite Elements. \emph{Computers and Structures 1991}; \textbf{40}:1043--1047.

\bibitem{Debongnie_1995} Debongnie~JF, Zhong~HG, Beckers P. Dual Analysis with General Boundary Conditions. \emph{Computer Methods in Applied Mechanics and Engineering 1995}; \textbf{122}:183--192.

\bibitem{Santos_2014} Santos~HAFA, Moitinho~de~Almeida~JP. A Family of Piola Kirchhoff Hybrid stress Finite Elements for Two-dimensional Linear Elasticity. \emph{Finite Elements in Analysis and Design 2014}; \textbf{85}:33--49.

\bibitem{Kempeneers_2010} Kempeneers~M, Debongnie~JF, Beckers P. Pure Equilibrium Tetrahedral Finite Elements for Global Error
Estimation by Dual Analysis. \emph{International Journal for Numerical Methods in Engineering 2010}; \textbf{81}:513--536.

\bibitem{Wang_2014} Wang~L, Zhong~H. A Traction-based Equilibrium Finite Element Free from Spurious Kinematic Modes for Linear Elasticity Problems. \emph{International Journal for Numerical Methods in Engineering 2014}; \textbf{99}:763--788.

\bibitem{BookAlmeidaMaunder} Almeida~JPM, Maunder~EAW. \emph{Equilibrium Finite Element Formulations}. John Wiley \& Sons Ltd. 2017.

\bibitem{AinsworthOden} \KRO{Ainsworth~M, Oden~TJ. \emph{A Posteriori Estimation in Finite Element Analysis}.
Wiley Interscience. 2000.}

\bibitem{LadevezeMaunder} \KRO{Ladev\`{e}ze~P, Maunder~EAW. A General Method for Recovering Equilibrating Element Tractions.
\emph{Computer Methods in Applied Mechanics and Engineering 1996}, \textbf{137}:111--151.}

\bibitem{ParesSantosDiez} \KRO{Par\'{e}s~N, Santos~H, D\'{\i}ez~P.
Guaranteed Energy Bounds for the Poisson Equation using a Flux-Free Approach: Solving the Local Problems in Subdomains.
\emph{International Journal for Numerical Methods in Engineering 2009}, \textbf{79}:1203--1244.}

\bibitem{LadevezeLeguillon} \KRO{Ladev\`{e}ze~P, Leguillon,~D.
Error Estimate Procedure in the Finite Element Method and Applications.
\emph{SIAM Journal on Numerical Analysis 1983}, \textbf{20}:485--509.}

\bibitem{BochevHyman} Bochev~BB, Hyman~JM. Principles of Mimetic Discretizations of Differential Operators. \emph{Compatible Spatial Discretizations, Eds.: D. Arnold, P. Bochev, R. Nicolaides and M. Shashkov, The IMA volumes in Mathematics and its Applications 2006}; \textbf{142}:89--119.

\bibitem{BonelleErn} Bonelle~J, Ern~A. Analysis of Compatible Discrete Operator Schemes for Elliptic Problems on Polyhedral Meshes. \emph{ESAIM Mathematical Modelling and Numerical Analysis 2014}; \textbf{48}(2):553--581.

\bibitem{Tonti} Tonti~E. Why Starting from Differential Equations for Computational Physics?. \emph{Journal of Computational Physics 2014}; \textbf{257}:1260--1290.

\bibitem{Kreeft_2011} Kreeft~J, Palha~A, Gerritsma~M. Mimetic Framework on Curvilinear Quadrilaterals of Arbitrary Order. \emph{Arxiv preprint 2011}.

\bibitem{Kreeft_stokes} Kreeft~J, Gerritsma~M. Mixed Mimetic Spectral Element Method for Stokes Flow: A Pointwise Divergence-free Solution. \emph{Journal of Computational Physics 2013}; \textbf{240}:284--309.

\bibitem{Douglas_2006} Arnold~DN, Falk~RS, Winther~R. Differential Complexes and Stability of Finite Element Methods II: The Elasticity Complex. \emph{Compatible Spatial Discretizations, Eds.: D. Arnold, P. Bochev, R. Nicolaides and M. Shashkov, The IMA volumes in Mathematics and its Applications 2006}; \textbf{142}:47--67.

\bibitem{Yavari_2008} Yavari~A. On Geometric Discretization of Elasticity. \emph{Journal of Mathematical Physics 2008}; \textbf{49}.

\bibitem{Angoshtari_2013} Angoshtari~A, Yavari~A. A Geometric Structure-preserving Discretization Scheme for Incompressible Linearized Elasticity. \emph{Computer Methods in Applied Mechanics and Engineering 2013}; \textbf{259}:130--153.

\bibitem{Yavari_2013} Yavari~A, Goriely~A. Nonlinear Elastic Inclusions in Isotropic Solids. \emph{Proceedings of the Royal Society A 2013}; \textbf{469}.

\bibitem{Angoshtari_2014} Angoshtari~A, Yavari~A. Differential Complexes in Continuum Mechanics. \emph{Archive for Rational Mechanics and Analysis 2014}; \textbf{216}:193--220.

\bibitem{Timoshenko_1982} Timoshenko~SP, Goodier~JN. \emph{Theory of Elasticity} (3rd~edn).
McGraw-Hill Inc. 1982.

\bibitem{Bertoti} Bertoti~E, Dual-mixed $p$ and $hp$ Finite elements for Elastic Membrane Problems, \emph{International Journal for Numerical Methods in Engineering 2002}, \textbf{53}: 3--29.

\bibitem{MEEVC} Palha A, Gerritsma MI. A Mass, Energy, Enstrophy and Vorticity Conserving (MEEVC) Mimetic Spectral Element Discretization for the 2D Incompressible Navier-Stokes Equations \emph{Journal of Computational Physics 2017}, \textbf{328}:200--220.

\bibitem{GerritsmaPhillips} Gerritsma MI, Phillips TN, Compatible Spectral Approximations for the Velocity-Pressure-Stress Formulation of the Stokes Problem, \emph{SIAM Journal of Scientific Computing 1999}, \textbf{20}:1530--1550.

\bibitem{Gerritsma_2011} Gerritsma~M. Edge Functions for Spectral Element Methods. \emph{Spectral and High Order Methods for Partial Differential Equations - Lecture Notes in Computational Science and Engineering 2010}; \textbf{76}:199--207.

\bibitem{Bossavit} Tarhasaari~T, Kettunen~L, Bossavit~A. Some Realizations of a Discrete Hodge Operator: A Reinterpretation of Finite Element Techniques. \emph{IEEE Transactions on Magnetics 1999}; \textbf{35}(3):1494--1497.

\bibitem{Hirani} Hirani~AN. Discrete Exterior Calculus. \emph{PhD Thesis, California Institute of Technology 2003}.

\bibitem{Lipnikova_2014} Lipnikov~K, Manzini~G, Shashkov~M. Mimetic Finite Difference Method. \emph{Journal of Computational Physics 2014}; \textbf{257}:1163--1227.

\bibitem{Canuto_2006} Canuto~C, Hussaini~M, Quarteroni~A, Zang~T. \emph{Spectral Methods, Fundamentals in Single Domains}
Springer. 2006.

\bibitem{Karniadakis_2005} Karniadakis~GE, Spencer~SJ. \emph{Spectral/hp Element Methods for Computational Fluid Dynamics} (2nd~edn).
Oxford Science Publications. 2005.

\bibitem{Reddy_2013} Reddy~JN. \emph{An Introduction to Continuum Mechanics} (2nd~edn).
Cambridge University Press. 2013.

\bibitem{Gerritsma_2012} Gerritsma~MI. An Introduction to a Compatible Spectral Discretization Method. \emph{Mechanics of Advanced Materials and Structures 2012}; \textbf{19}(3):48--67.

\bibitem{Gordon_1973} Gordon~WJ, Hall~CA. Construction of Curvilinear Coordinate Systems and Applications to Mesh Generation. \emph{International Journal for Numerical Methods in Engineering 1973}; \textbf{7}:461--477.

\bibitem{bookPian} Pian~THH, Wu~C-C. \emph{Hybrid and Incompatible Finite Element Methods}, Chapman \& Hall/CRC, 2006.
    
\bibitem{BrezziFortin} Brezzi~F, Fortin~M.
    \emph{Mixed and Hybrid Finite Element Methods}, Springer Series in Computational Mathematics 15, 1991.   
    
\bibitem{BoffiBrezziFortin} Boffi~D, Brezzi~F, Fortin~M.
    \emph{Mixed Finite Element Methods and Applications}, Springer Series in Computational Mathematics 44, 2013.
    
\bibitem{Jain} \KRO{Jain~V, Zhang~Y, Palha~A, Gerritsma~M.
    Construction and Application of Algebraic Dual Polynomial Representations for Finite Element Methods. \emph{	arXiv:1712.09472 [math.NA] 2017}.}
\end{thebibliography}
\end{document}